\documentclass[12pt]{amsart}
\usepackage{amsthm,amsmath,amssymb,amscd,graphics,enumerate}

\newcommand{\marg}[1]{\footnote{#1}\marginpar{\tiny\thefootnote}}
\newcommand{\thisdate}{\today}
{
   \newtheorem{theorem}[subsubsection]{Theorem}
   \newtheorem{proposition}[subsubsection]{Proposition}     
   \newtheorem{lemma}[subsubsection]{Lemma}
   \newtheorem*{claim}{Claim}
   \newtheorem{corollary}[subsubsection]{Corollary}

}
{\theoremstyle{definition}

   \newtheorem{definition}[subsubsection]{Definition}
   \newtheorem{remark}[subsubsection]{Remark}
}
\newcommand{\RR}{{\mathbb{R}}}
\newcommand{\bbS}{{\mathbb{S}}}

\newcommand{\CC}{{\mathbb{C}}}

\newcommand{\ZZ}{{\mathbb{Z}}}

\newcommand{\GG}{{\mathbb{G}}}
\newcommand{\LL}{{\mathbb{L}}}

\newcommand{\bB}{{\mathbf{B}}}

\newcommand{\bK}{{\mathbf{K}}}
\newcommand{\bM}{{\mathbf{M}}}
\newcommand{\bR}{{\mathbf{R}}}
\newcommand{\bX}{{\mathbf{X}}}

\newcommand{\bmu}{{\boldsymbol{\mu}}}
\newcommand{\cA}{{\mathcal A}}
\newcommand{\cB}{{\mathcal B}}
\newcommand{\cC}{{\mathcal C}}

\newcommand{\cE}{{\mathcal E}}

\newcommand{\cH}{{\mathcal H}}

\newcommand{\cK}{{\mathcal K}}

\newcommand{\cM}{{\mathcal M}}
\newcommand{\cN}{{\mathcal N}}
\newcommand{\cO}{{\mathcal O}}

\newcommand{\cS}{{\mathcal S}}
\newcommand{\cT}{{\mathcal T}}

\newcommand{\cX}{{\mathcal X}}
\newcommand{\cY}{{\mathcal Y}}

\newcommand{\pre}{\operatorname{pre}}

\newcommand{\Spec}{\operatorname{Spec}}

\newcommand{\Isom}{\operatorname{Isom}}
\newcommand{\underisom}{{\underline{\operatorname{Isom}}}}
\newcommand{\bIsom}{\mathbf{Isom}}

\newcommand{\Hom}{{\operatorname{Hom}}}
\newcommand{\cHom}{{{\cH}om}}
\newcommand{\underhom}{{\underline{\operatorname{Hom}}}}

\newcommand{\Ext}{{\operatorname{Ext}}}
\newcommand{\cExt}{{{\cE}xt}}

\newcommand{\Ker}{{\operatorname{Ker}}}
\newcommand{\Aut}{{\operatorname{Aut}}}
\newcommand{\underaut}{{\underline{\operatorname{Aut}}}}
\newcommand{\Out}{{\operatorname{Out}}}

\newcommand{\Def}{{\operatorname{Def}}}
\newcommand{\sh}{{\operatorname{sh}}}

\newcommand{\opH}{{\operatorname{H}}}

\newcommand{\cyl}{\operatorname{cyl}}

\newcommand{\lrar}{\longrightarrow}
\newcommand{\dar}{\downarrow}
\newcommand{\down}{\downarrow}

\newcommand{\ocM}{\overline{{\mathcal M}}}
\newcommand{\obM}{\overline{{\mathbf M}}}

\newcommand{\oC}{{\overline{C}}}
\newcommand{\oDelta}{{\overline{\Delta}}}

\newcommand{\sing}{_{\operatorname{sing}}}
\newcommand{\bal}{{\operatorname{bal}}}

\newcommand{\rig}{{\operatorname{rig}}}
\newcommand{\tei}{{\operatorname{tei}}}

\newcommand{\Adm}{{\mathcal A}dm}

\newcommand{\Zm}{\ZZ/m\ZZ}
\newcommand{\level}{(\Zm)^{2g}}
\newcommand{\double}{\genfrac..{0pt}1
{\raise -1pt\hbox{$\scriptstyle\longrightarrow$}}{\raise 3pt\hbox
{$\scriptstyle\longrightarrow$}}} 

\newcommand{\prym}[2]{{\genfrac[]{0pt}1{_{^{#1}}}%
{^{_{#2}}}}}   

\newcommand{\KO}[4]{{\cK_{#1,#2}(#3,#4)}}
\newcommand{\KOB}[4]{{\cK^{\bal}_{#1,#2}(#3,#4)}}
\newcommand{\ko}[4]{{\bK_{#1,#2}(#3,#4)}}

\newcommand{\TSM}{\KO{g}{n}{\cM}{d}}
\newcommand{\tsm}{\ko{g}{n}{\cM}{d}}
\newcommand{\SM}{\KO{g}{n}{\bM}{d}}
\newcommand{\sm}{\ko{g}{n}{\bM}{d}}
\newcommand{\TSMB}{\KOB{g}{n}{\cM}{d}}

\newcommand{\TC}[3]{{\mathcal{B}^{#2}_{#1}(#3)}}

\newcommand{\TCbal}[1]{\TC{g,n}{\operatorname{bal}}{#1}}

\newcommand{\TCtei}[1]{\TC{g}{\operatorname{tei}}{#1}}

\def\gen{_{\rm gen}}
\def\rest{|_}

\def\tototi{\mathbin{\mathop{\otimes}\limits^{\raise-1pt\hbox
{$\scriptscriptstyle {\rm L}$}}}}

\def\indlim{\mathop{\vrule width0pt height7pt depth
4pt\smash{\lim\limits_{\raise 1pt\hbox to 14.5pt
{\rightarrowfill}}}}}
\def\projlim{\mathop{\vrule width0pt height7pt depth
4pt\smash{\lim\limits_{\raise 1pt\hbox to 14.5pt
{\leftarrowfill}}}}}

\def\mapright#1{\mathrel {\smash{\mathop{\longrightarrow}
\limits^{#1}}}}

\def\down{\big\downarrow}

\def\mapdown#1{\down\rlap{$\vcenter{\hbox
{$\scriptstyle#1$}}$}}
\def\into{\hookrightarrow}
\begin{document}
\title[Twisted bundles]{Twisted bundles and admissible covers}
\author[D. Abramovich]{Dan Abramovich}
\thanks{D.A. Partially supported by NSF grant DMS-9700520 and by an Alfred
	P. Sloan research fellowship}  
\address{Department of Mathematics\\ Boston University\\ 111 Cummington
	 Street\\ Boston, MA 02215\\ USA} 
\email{abrmovic@math.bu.edu}
\author[A. Corti]{Alessio Corti}
\address{DPMMS, CMS\\ 
University of Cambridge\\ 
Wilberforce Road\\
Cambridge CB3 0WB\\ 
U.K.} 
\email{a.corti@dpmms.cam.ac.uk}
\author[A. Vistoli]{Angelo Vistoli}
\thanks{A.V.  partially supported by the University of
	Bologna, funds for selected research topics.}
\address{Dipartimento di Matematica\\ Universit\`a di Bologna\\Piazza di Porta
	 San Donato 5\\ 40127 Bologna\\ Italy}
\email{vistoli@dm.unibo.it}
\date{\thisdate}

\maketitle

\setcounter{tocdepth}{1}
\tableofcontents
\section{Introduction}

The purpose of  this paper is twofold. First, we discuss and prove
results on twisted covers announced without proofs in \cite{modfam},
section 3 (with slightly modified  notation). Second,   
we continue with some new results about nonabelian level structures.  
 Some results related to ideas in this paper
were also discovered independently by F. Wewers
\cite{Wewers}. 

For the sake of motivation, we start with the problem of smooth moduli
spaces of stable curves with level structures.

\subsection{Moduli of curves - course and fine}
The moduli space of smooth curves $\bM_g$ and its natural
compactificaion by the moduli space of stable curves $\obM_g$, are
among the most illustrious successes of twentieth century
mathematics. Yet they have a somewhat unpleasant feature - they are
in general singular, where they really shouldn't be - the deformation
spaces of stable curves are all smooth, yet automorphisms prevent the
space from being a fine moduli space, and often force the
coarse moduli spaces to be singular. Nowadays one knows that this is
in some sense an ``optical illusion'' - one should really work with
the corresponding moduli stack, which is always nonsingular. Yet it is
a bit dissatisfying to require the use of a specialized tool-kit such
as algebraic stacks to see the smoothness in a moduli problem which
should have been visibly smooth for the bare eyes.

A satisfactory solution for the ``open'' moduli space was proposed by
Mumford -  the space $\bM_g$ admist a finite Galois cover $\bM_g^{(m)}$, the
moduli space of curves with level-$m$ structure, and as soon as $m\geq
3$ this is a smooth, fine moduli space.

The search for a similar solution for  $\obM_g$ has taken several
turns through the years. In \cite{Mumford} Mumford used the normalization
$\obM_g^{(m)}$ of  $\obM_g$ in  $\bM_g^{(m)}$:
$$ \begin{array}{ccc}\bM_g^{(m)} & \subset & \obM_g^{(m)}\\
\dar & & \dar\\
                      \bM_g      & \subset & \obM_g
\end{array}$$
as a finite covering of
$\obM_g$ which carries a ``tautological family'' of stable curves. This
covering is in general singular, it is not a fine moduli space of a
natural moduli problem, yet the fact that its singularities are
Cohen-Macaulay was useful for intersection theory.

Looijenga \cite{Looijenga}, and soon after Pikaart and De Jong
\cite{Pikaart-DeJong} (see also \cite{Boggi-Pikaart}), used instead
the normalizations 
$\overline{{}_G\bM_g}$ of $\obM_g$ in the moduli spaces of smooth curves with
Teichm\"uller level structures ${{}_G\bM_g}$. They showed that, with
careful appropriate  choices of finite groups $G$, these spaces are
smooth Galois   
covers of $\obM_g$. Yet they did not write down a simple description
of these spaces as fine moduli spaces of appropriate moduli problems. 

It seems that the main reason a fine--moduli--space interpretation of
$\overline{{}_G\bM_g}$ was not given in \cite{Looijenga},
\cite{Pikaart-DeJong} is that this would require working with certain
objects up to outer automorphisms. We describe this situation in this paper in
terms of a process of {\em rigidification} of a stack. It follows that
the spaces studied by Looijenga and Pikaart--De Jong are indeed fine
moduli spaces of certain Teichm\"uller  structures obtained by
rigidification. This solves the problem, but one would prefer a
solution which  avoids
the process of rigidification altogether. 

In Section \ref{Sec:fine-covers} we introduce a moduli space which entirely
circumvents the need for rigidification. For any $g\geq 2$ we
give a finite group $G$ 
such that the moduli space of connected admissible $G$ covers of genus
$g$ is a smooth, fine moduli space, which is a Galois cover of $\obM_g$.
The proofs rely on methods introduced in \cite{Looijenga} and
\cite{Pikaart-DeJong}, and on our theory of  twisted
$G$-covers, developed in the first few sections of the paper and
summarized in \ref{Sec:summary} below.

\subsection{Algebra-to-analysis translation
table}\label{Sec:algebra-to-analysis}\hfill 

In this paper we systematically use the language of stacks.
However,
the results in this paper may be of interest for people studying moduli
spaces from the point of view of differential geometry and
analysis, some of whom may prefer the language of differential
orbifolds. These people should be able to get by using the following
rough translation table.\\

{\small
\begin{tabular}{|l|l|}\hline
Schemes over a base scheme $\bbS$ & Analytic spaces (over $\Spec \CC$)
\\ \hline
Smooth algebraic curve $C \to S$& Family of Riemnann surfaces $C\to S$
\\ \hline
Stable algebraic curve $C \to S$& Family of stable Riemnann surfaces $C\to S$
\\ \hline
Geometric point $s:\Spec \Omega \to X$ &
\parbox{3in}{\begin{enumerate} \item a point of $X$ (if $\Omega 
= \CC$)  \item  a general point of an
irreducible closed subset (general case)\end{enumerate}}
\\ \hline
An \'etale neighborhood of a point &  A neighborhood in the Euclidean topology
\\ \hline
\parbox{3in}{The strict henselization $X^\sh$ at a geometric point} &
\parbox{3in}{A {\em small,  contractible} 
neighborhood of a point of $X$ -- or -- the germ of $X$ at a point}
\\ \hline
A category which is an algebraic stack $\cX$ & \parbox{3in}{A reasonably good
moduli problem of families of geometric objects } 
\\ \hline
A (tame) Deligne--Mumford stack $\cX$ & \parbox{3in}{An ``orbifold'' with
finite stabilizers,
possibly with the local groups action not effective}
\\ \hline
The stack theoretic quotient $\cX = [U/G]$ & The space $U/G$ along
with its orbifold structure
\\ \hline
A moduli Deligne--Mumford stack $\cX$ & \parbox{3in}{A moduli orbifold
representing 
the moduli 
problem $\cX$.}
\\ \hline
An object in $\cX(S)$ & \parbox{3in}{A family $Y \to S$ of object in
the moduli problem}  
\\ \hline
A 1-morphism of Deligne--Mumford stacks & \parbox{3in}{A ``good map''
of orbifolds in 
the sense of Chen-Ruan, without the assumption of injectivity of stabilizers}
\\ \hline
A representable 1-morphism & \parbox{3in}{A ``good map''
of orbifolds in 
the sense of Chen-Ruan}
\\ \hline
The coarse moduli space $\bX$ of $\cX$ & The analytic space $\bX$
underlying an orbifold $\cX$ 
\\ \hline
\end{tabular}}

It is time to introduce the main concepts studied in this paper.

\subsection{Summary of the paper}\label{Sec:summary}

We fix a noetherian base scheme $\bbS$. The analytically-inclined
reader may very well assume that $\bbS = \Spec \CC$ and work
complex-analytically, perhaps using the translation table provided
in section \ref{Sec:algebra-to-analysis} above. 

\subsubsection{Twisted stable maps} Kontsevich's stack
$\ocM_{g,n}(X,\beta)$ of $n$-pointed stable 
maps of genus $g$ into a projective scheme $X$ with homology class
$\beta$ (see \cite{Kontsevich}) have served as an extremely useful tool
in enumerative 
geometry (see e.g. \cite{Fulton-Pandharipande}) and as a construction
technique (see e.g., \cite{A-O:Hurwitz}). For
similar reasons it is of interest to replace $X$ by a Deligne--Mumford stack.
In the paper \cite{stablemaps}, a proper stack $\TSM$ of {\em twisted stable
 maps}  was constructed for each tame Deligne-Mumford stack $\cM$
admitting a projective coarse moduli scheme    $\bM$. 
 As it turns out, the stack of ``usual'' stable maps into a
 stack $\cM$ may fail to be  proper -  an example relevant to this
paper is given in \cite{stablemaps}, section 1.3. The main point    
of the paper \cite{stablemaps} was, that in order to have a {\em proper}
stack of maps into the stack $\cM$, it is natural to allow the curves
 at the source of the 
 map to acquire orbispace structure at the nodes (and to
 have a complete picture, also along
 the markings). A particularly important open-and-closed substack of
 $\TSM$ is the stack $\TSMB$ of {\em balanced} twisted stable maps, where the
 underlying twisted curves are smoothable.

Here we study the constructions of
\cite{stablemaps} in the case where $\cM$ is the classifying stack
$\cB G$ of a finite group (or, more generally, a finite \'etale group
scheme) $G$. (The  assumption that $\cB G$ be tame translates to the
requirement that for each field $k$ and point $p$ in $\bbS(k)$,  the
degree of $G_p$ is invertible in $k$.) 

\subsubsection{Twisted stable maps to $\cB G$ via $G$-covers}
There are several reasons why we find it interesting to
pursue this special case. On the one hand, restricting to this special
case allows us to give   fairly explicit descriptions of the stacks of
twisted stable maps, notably in terms of certain ramified Galois
$G$-covers of 
``classical'' stable pointed curves, with no reference to orbifold
curves of any kind. Constructing the stack of these $G$-covers
directly is  a standard procedure in moduli theory, and thus
one can completely circumvent the delicate steps that were used in
\cite{stablemaps} to construct the stack of twisted stable maps.
\marg{In the final version an appendix on such construction will be added}

One
hopes  that one might be able to use ideas of this paper to
significantly simplify the constructions of \cite{stablemaps} in
general. There is also some hope that further study of this case may shed
light on possible constructions for a non-tame target stack $\cM$, see
\cite{A-O:Hurwitz}, \cite{Bouw-Wewers}.
Discussion of twisted stable maps in this special case of $\cB G$ also
serves as 
an opportunity  to study some properties of twisted curves
and twisted stable maps which were not addressed in \cite{stablemaps}. 

\subsubsection{Hurwitz stacks  via twisted stable maps} 
On the other hand, having the tool-kit of twisted stable maps
available allows us to shed new light on a time-honored topic in
algebraic geometry: the stacks of {\em balanced} twisted $G$-covers
are easily interpreted as a certain compactification of a Hurwitz-type stack
of Galois covers of curves. 

First, it is an easy exercise in deformation theory to show
that this stack of balanced twisted covers is always nonsingular and
of dimension $3g-3+n$, 
giving a ``finite'' (though in general not representable) flat cover
of $\ocM_{g,n}$. 

Restricting to the case where
$G$ is the symmetric group, and using the usual correspondence between
\'etale covers of degree $n$ and $\cS_n$-covers, we get a
compactification of a ``usual'' Hurwitz-type  stack (without group
actions). We 
show that this is no other than the {\em normalization} of the
Harris--Mumford stack of admissible covers (with arbitrary ramification
type). Giving a moduli 
interpretation to this normalization has been a desirable goal since
the appearance of \cite{Harris-Mumford}, since the singularities of
the stack of admissible covers are not natural for  most
applications.

\subsubsection{Rigidification and Teichm\"uller structures}
We proceed to consider the open-and-closed substack of
$\TCbal{G}$ consisting of {\em connected} balanced twisted $G$-covers of
unpointed stable curves of genus $g$. A
connected $G$-bundle over a smooth curve corresponds to a
Teichm\"uller level structure in the sense of \cite{Deligne-Mumford},
Section 5, however the center  $Z(G)$ of $G$ acts on the
$G$-bundle fixing 
the level structure. We describe a fairly general procedure of
removing an \'etale group action from the center of stabilizers in a
stack, which we call {\em rigidification}. Thus the rigidification of
$\TCbal{G}$ obtained by removing the center $Z(G)$ is a smooth, proper
Deligne--Mumford stack $\TCtei{G}$, which serves as a natural
compactification of the stack ${}_{G}\cM_g$ of Teichm\"uller level
structures on  
smooth curves. We call this stack {\em the stack of twisted
Teichm\"uller structures}. This stack deserves the notation
$\overline{{}_G\cM_g}$, but unfortunately this notation has been wrongfully
used for another stack by other authors (Deligne, Mumford, Pikaart, De
Jong \ldots).

\subsubsection{A good compactification of $\bM_g^{(m)}$}
We find that two special situations are of particular interest. First, if $G =
(\ZZ/m\ZZ)^{2g}$, then a Teichm\"uller $G$-structure on a smooth
curve corresponds to a level-$m$ structure in the usual sense, and
$\TCtei{G}$ is a natural compactification of Mumford's space of curves
with level-$m$  structure. We give a detailed description of the type
of objects that appear in the boundary: the precise ``twisting'' of
the underlying twisted curves, and an interpretation of the twisted
level-$m$ structure in terms of a trivialization of the \'etale
cohomology with values in $\ZZ/m\ZZ$ of the twisted curve.
If the structure sheaf of the
base scheme contains the $m$-th roots of 1, then one can also define a
symplectic structure on the \'etale cohomology group.\marg{In the final
version, a
discussion of symplectic structure wil be added} 

Mumford considered a different compactification - the normalization
$\obM_g^{(m)}$ of
$\obM_g$ in the fuction field  of the space $\bM_g^{(m)}$  of smooth
curves with level 
structures. Unlike $\TCtei{G}$, this space is singular, a fact which
here we see 
resulting from the fact that a twisted level-$m$ structure on a
singular curve may have automorphisms. Indeed $\obM_g^{(m)}$ is the
coarse moduli space of $\TCtei{G}$. 
Our setup allows us to show, as a fairly easy consequence of Serre's
lemma, the well known fact that 
Mumford's 
compactification does carry a ``tautological family'' of stable curves 
(albeit without a level structure), i.e. a morphism $\obM_g^{(m)}\to \ocM_g$. 

\subsubsection{A projective fine moduli space of curves with level structure}
The second special situation we consider is that of ``rigid''
nonabelian level structures, in particular those discussed in the papers
\cite{Looijenga} of E.\ Looijenga  and
\cite{Pikaart-DeJong} of M.\ Pikaart and A.J.\ de Jong. In
\cite{Looijenga} the group in question is the structure group of a
$\prym{2m}{2}$-Prym level structure. In
\cite{Pikaart-DeJong}
the group $G$ is 
the maximal nilpotent quotient of exponent $n$ and nilpotence order $k$ of the
fundamental group of a Riamann surface of genus $g$. Consider the
normalization $\overline{{}_G\bM_g}$ of $\obM_g$ in the space of smooth
curves with Teichm\"uller level-$G$ structure ${}_G\bM_g$. Looijenga
and Pikaart and De 
Jong show in their respective cases that, for suitable values of the
parameters $m$, respectively $(n,k)$, the space  $\overline{{}_G\bM_g}$    
is a {\em nonsingular} finite cover of $\obM_g$. In this case we are
able to show that this space {\em coincides} with our stack
$\TCtei{G}$; thus its non-singularity can be interpreted in terms of
the fact that the automorphism group of a twisted Teichm\"uller
$G$-structure is always trivial. 

We improve on these results by introducing, for each genus $g$, a
group $G$, such that moreover the
automorphism group of every connected admissible $G$-cover is
trivial (in particular $G$ has trivial center). The group $G$ can be
quickly described as follows: if $p_1$ 
and $p_2$ are two distinct primes, let $G_1$ be the structure group of
a $\prym{p_1}{p_2}$-Prym level structure, and let $G_2$ be the
structure group of 
a $\prym{p_2}{p_1}$-Prym level structure. These groups admit natural
homomorphisms to $(\ZZ/p_1p_2\ZZ)^{2g}$, and $G$ is the fibered product.

With this $G$, our stack $\TCtei{G}$ is then a fine moduli space of
connected, admissible $G$ covers, and the rigidification process is
avoided.

\subsection{Summary of notation}

\begin{itemize}
\item $\bbS$ ---  the base scheme.
\item $\cC$ --- a twisted curve. 
\item $\pi:\cC \to C$ --- the morphism to the coarse moduli scheme.
\item $C\gen$ --- the generic logus of $C$ (or $\cC$).
\item $\Sigma_i^\cC, \Sigma_i^C$ --- the $i$-th marking of $\cC$ and
$C$, the union of which is denoted $\Sigma^\cC$, respectively, $\Sigma^C$.
\item $G$ --- a finite group.
\item $\cB G$ --- the classifying stack of $G$.
\item $\TC{g,n}{}{G}$ --- the stack of twisted $G$-covers of
$n$-pointed curves of genus $g$.
\item $\TC{g}{}{G}$ --- the stack of twisted $G$-covers of (unpointed)
curves of genus $g$. 
\item $\TCbal{G}$ --- the stack of {\em balanced} twisted $G$-covers.
\item $\TCtei{G}$ --- the stack of twisted Teichm\"uller
$G$-structures (of unpointed curves) .
\item $\Adm_{g,n,d}$ --- the stack of generalized Harris--Mumford
admissible covers of  degree $d$ over $n$-pointed curves of genus
$g$. 
\item $\Adm_{g,n}(G)$ --- the stack of admissible $G$-covers.
\end{itemize} 

\subsection{What we mean by ``the local picture''} We often need to
``identify'' a scheme, or a diagram connecting several schemes,
locally in terms of explicit equations. To avoid repeated mention of
\'etale localizations or strict henselizations,  we make the following
agreement: we say that ``the local picture of $X$ at a geometric point $p$ is
the same as $U$ (at point $q$)'' if  the   germ of $X$ at $p$ is
isomorphic to the germ of 
$U$ at $q$, in other words: 
\begin{enumerate} 
\item in the algebraic situation: there is an isomorphism between the
strict henselization $X^\sh$ of $X$ at $p$ and the strict
henselization  $U^\sh$ of $U$ at $q$. Often, however,  it is enough to
assume that there is an \'etale neighborhood $X'$ of $p$ and an
isomorphism $X' \tilde{\to} U'$ with an \'etale neighborhood $U'$ of $q$.
\item In the analytic situation: there is a small contractible
neighborhood $X'$ of $p$ and an
isomorphism $X' \tilde{\to} U'$ with a neighborhood $U'$ of $q$.
\end{enumerate}

\subsection{Acknowledgements} We thank Johan de Jong for extremely
helpful ideas 
regarding non-abelian level structures.

\section{Terminology}

In this section we recall basic facts about twisted stable maps, and
introduce the notion of twisted $G$-covers.

\subsection{Twisted curves and twisted stable maps}

We follow the setup in \cite{stablemaps} for defining twisted stable
maps. The basic object underlying a 
twisted stable map is a {\em twisted curve}, i.e., a pointed nodal
curve $C$ along 
with a Deligne--Mumford stack structure $\cC$ at its nodes and
markings. 

The local picture of a twisted curve at a geometric point $p$ can be
explicitly described as follows: 

\subsubsection{At a marking}\label{Sec:at-marking}
The local picture of $\cC \to S$ is the same as $[U/\mu_r]\to T$, where 
\begin{enumerate}
\item $T  =    \Spec A$, 
\item  $U  =  \Spec A[z]$, 
and 
\item the action of $\mu_r$ is given by $z \mapsto \zeta_r \cdot z$.
\end{enumerate}

\subsubsection{At a node}\label{Sec:at-node} If $p$ lies over a node,
the local 
picture of $\cC \to S$ is the same as $[U/\mu_r]\to T$, where  
\begin{enumerate}
\item $T  =    \Spec A$,
\item  $U  =   \Spec A[z,w]/(zw-t)$
for some $t\in A$, and 
\item\label{Item:node-action} the action of $\mu_r$ is given by $(z,w)
\mapsto (\zeta_r z, \zeta_r^aw)$, for some $a\in (\ZZ/r\ZZ)^\times$. 
\end{enumerate}
This description 
is implicit in \cite{stablemaps}, Proposition 3.2.3.

\subsubsection{Balanced actions, balanced curves} Note that, unless $a
\equiv -1 \mod r$ in \ref{Sec:at-node}(\ref{Item:node-action}) above,
the element $t$  
must vanish, and the node cannot be smoothly deformed; therefore the
``locally smoothable'' case $a\equiv -1 \mod r$ deserves special
attention. An action with 
$a \equiv -1 \mod r$ is called {\em  
balanced}, a node presented via a balanced action is called balanced,
and a twisted curve is said to be balanced when all its nodes are
balanced.


Here is a formal definition of twisted curves:

\begin{definition}
A {\em twisted nodal $n$-pointed curve over a scheme $S$} is a  diagram
$$\begin{array}{ccc} \Sigma^{\cC} & \subset & \cC  \\  &\searrow & \dar  \\ &
 & C   
  \\ &  & \dar  \\ &&S  
\end{array}$$
where
\begin{enumerate}
\item $\cC$ is a tame Deligne-Mumford stack, proper over $S$, which
\'etale locally is a nodal curve 
over $S$;
\item $\Sigma^\cC = \cup_{i=1}^n \Sigma_i^\cC$, where $\Sigma_i^{\cC}
\subset  \cC$ are disjoint closed substacks in the 
smooth locus 
of $\cC \to S$;  
\item $\Sigma_i^{\cC} \to S$ are \'etale gerbes; 
\item the morphism $\cC \to C$ exhibit $C$ as the coarse moduli scheme of 
$\cC$; and
\item $\cC \to C$ is an isomorphism over $C\gen$.
\end{enumerate}
\end{definition}

The notation $C\gen$ in the definition above stands for the {\em
generic locus}, namely the complement of the nodes and markings on
$C$. We denote the image marking of $\Sigma_i^{\cC}$ in $C$ by
$\Sigma_i^C$.
The stabilizer of a geometric point $p$ over a node or marking will be
denoted $\Gamma_p$.

\subsubsection{Twisted stable maps}
Fix a proper deligne-Mumford stack $\cM$ over $\bbS$,
admitting a projective coarse moduli scheme $\bM$, on which we fix a very ample
invertible sheaf. For arithmetic situations we need a bit more: as in
\cite{stablemaps}, we assume that $\cM$ is tame, namely the order of
the automorphism group of a geometric object of $\cM$ is prime to the
residue characteristic of the field over which the object is defined.

We can now recall the definition of a twisted stable map:

\begin{definition}\label{Def:twisted-stable-map} A {\em twisted stable
$n$-pointed map of genus $g$ and degree $d$ over $S$}
$$(\cC \to S, \Sigma^{\cC}\subset \cC, f\colon  \cC \to \cM)$$   consists
of a commutative diagram
$$\begin{array}{ccc} \cC &\to& \cM \\
		\dar & & \dar \\
		C & \to & \bM \\
\dar &&\\ S &&
\end{array}$$
 along with a closed substack $\Sigma^{\cC}\subset \cC$,
satisfying:
\begin{enumerate}
\item $\cC \to C \to S$ along with $\Sigma^{\cC}$  is a twisted nodal
$n$-pointed curve  over $S$; 
\item $(C\to S, \Sigma^C, f\colon C \to \bM)$ is a stable $n$-pointed map
of degree $d$; and   
\item the morphism $\cC \to \cM$ is representable.
\end{enumerate}
\end{definition}

While the first two conditions are ``classical'', the third is
truly stack-theoretic: it means that the stabilizer of a geometric point in
$\cC$ injects into the stabilizer of the image point in $\cM$. It can
be viewed as a stability condition - it is 
there to guarantee that the moduli problem be separated. In effect, we
allow $\cC$ to be just  twisted enough to afford a morphism to $\cM$.  

Twisted curves (and thus twisted stable maps) naturally form a
2-category, where the 1-morphisms are given by fiber diagrams. In
\cite{stablemaps} it was shown that the automorphism group of every
1-morphism is trivial, therefore this 2-category is equivalent to the
associated category, obtained by replacing 1-morphisms by their
2-isomorphism classes. The category of twisted stable maps thus obtained
is denoted $\TSM$. The main theorem of \cite{stablemaps} is:

\begin{theorem}\label{Th:stable-maps}
\begin{enumerate}
\item 
The category $\TSM$ is a proper algebraic stack with finite diagonal. 
\item The coarse moduli space
$\tsm$ of $\TSM$ is projective. 
\item There is a commutative diagram 
$$\begin{array}{ccc} \TSM & \to & \SM \\ 
\dar & & \dar \\ \tsm & \to & \sm
\end{array}
$$
where the top arrow is proper, quasifinite, relatively of
Deligne--Mumford type and tame, and the bottom arrow is finite. In particular,
if $\SM$ is a Deligne--Mumford stack, then so is
$\TSM$.
\end{enumerate} 
\end{theorem}

\subsubsection{Twisted objects} 
In the paper \cite{stablemaps}, a
considerable 
effort was made to give 
an explicit realization of $\TSM$ as a category of  twisted objects over
``usual'' curves, given using local charts.
  
 Briefly,
a twisted stable map gives rise to an object $\xi$ of
$\cM(C\gen)$. Also, \'etale locally we can present $\cC$ around a marking or a
node by  $[U/\Gamma]$, where $U$ is a (non-proper) marked
curve and 
$\Gamma$ is a finite group whose action on $U$  is free on
$U\gen$.  Over $U$ 
we have a $\Gamma$-equivariant object  $\eta \in \cM(U)$. The data
$(U, \eta,\Gamma)$ is called a {\em chart}.

 It was shown that
the collection of such data $(U, \eta,\Gamma)$ are compatible charts
in an atlas for a {\em twisted $\cM$-valued object} over
$C$. Further a moduli category of twisted $\cM$-valued
objects was defined. It was also shown that there is a base-preserving
equivalence of categories between $\TSM$ and the stack of twisted objects. The
category of twisted objects is closely related to the moduli problem
described by Chen and Ruan in \cite{Chen-Ruan}.

A slightly different, and somewhat simpler, realization is given below
in this
paper in the case where $\cM=\cB G$.

\subsection{Twisted $G$-covers}\label{Sec:twisted-G-covers}

We now introduce {\em twisted $G$-covers,} the main objects of our
paper, in terms of twisted stable maps. The reader who finds this
difficult to picture is encouraged to read Section
\ref{Sec:admissible-G-covers}, where  a concrete realization is given.

Let $G$ be a finite group, or, more generally, a finite \'etale
group scheme of constant degree. We assume that the degree is
prime to all residue characteristics in the base scheme
$\bbS$. 

Consider the classifying stack $\cB G$ of $G$. We recall that
$\cB G$ is a category whose objects
 over a scheme $T$ are principal $G_T$ bundles $P \to T$,
and morphisms are $G$-equivariant fiber diagrams of such principal
bundles:
$$ \begin{array}{ccc} 
P_1 & \to & P_2 \\
 \dar & & \dar \\   
T_1 & \to & T_2.
\end{array}$$

{\em We stick with the tradition that $G$ acts on a principal bundle $P \to
S$ on the right.} On the other hand, we write $G$-automorphisms of $P
\to S$ on the left:  if $S$ is connected this allows us to identify
these $G$-automorphisms as elements of $G$ rather than $G^{opp}$.

We have a presentation $\cB G = [\bbS /G]$, with $G$ acting trivially
on $\bbS$. Also useful is the
fact that 
the coarse moduli space of $\cB G$ is simply $\bbS$.

 Define a {\em  twisted $G$-cover of an $n$-pointed  curve of 
genus $g$} to be an object of the stack of twisted stable maps  $\KO{g}{n}{\cB
G}{0}$. According to Theorem \ref{Th:stable-maps} quoted above, this
stack is a proper stack whose moduli space is projective over
$\bbS$. Also,  since 
the coarse moduli scheme of $\cM = \cB G$ is the base scheme 
$\bM = \bbS$, the stack $\SM = \ocM_{g,n}$ is a Deligne--Mumford stack,
and therefore $\KO{g}{n}{\cB
G}{0}$ is also a Deligne Mumford stack. From
now on we  
will denote it by $\TC{g,n}{}{G}$. We will also denote by $\TCbal{G}$
the open 
and closed substack $\KOB{g}{n}{\cB G}{0}$ consisting of {\em
balanced} 
stable maps. 

In accordance to the name ``twisted $G$-cover'', an object of
$\TC{g,n}{}{G}$ over a scheme $S$ will be represented by 
the associated principal bundle $P \to \cC$, where $\cC\to S$ is the
underlying twisted curve. To avoid confusion, we will refer to  the morphism
 $\cC \to \cB G$ as the twisted stable map associated to  the twisted
$G$-cover $P \to \cC$.

We  discuss the structure of $P \to \cC$ in some detail later in
this paper. One useful fact we cite right away is the following:

\begin{lemma}\label{Lem:P-representable} $P\to S$ is a projective nodal curve.
\end{lemma}

\proof Note that $\cB G = [\bbS / G]$. The morphism $\bbS \to \cB G$,
the universal principal $G$ bundle, is clearly \'etale and finite.
It is also representable since $\bbS$ is. Also the morphism $\cC \to \cB G$ is
representable being a twisted stable map. Therefore $P = \cC
\times_{\cB G} \bbS$ is representable. For the same reason $P \to \cC$
is finite. Since $\cC \to C$ is quasi-finite and proper, so is the
composition $P \to C$. Since $C\to S$ is projective, it follows that
$P \to S$ is 
projective. Finally $P\to S$ is a nodal curve since it is \'etale over
the nodal twisted curve $\cC$. \qed

\section{Deformation theory of  twisted covers}\label{Sec:deformation-theory}

We show that the stack of twisted $G$ covers is unobstructed,
therefore smooth, and we calculate its diemsnsion. We note that, at
least in the balanced case, this can be shown  on the level of
Galois admissible covers, see \cite{Wewers}. A somewhat less detailed
version of this argument can be found in \cite{A-Jarvis}.

\begin{theorem} The stack  $\TC{g,n}{}{G}$ is
smooth. Its dimension
at a given twisted $G$-cover is $3g+n-3-u$, where $u$ is 
the number of nodal points at which the bundle in not 
balanced. 
\end{theorem}
 
{\bf Proof.} The cotangent complex $\LL_{\cB G/\bbS}$ of $\cB G$ is 
trivial, therefore deformations and obstructions  
of a twisted stable map $\cC \to \cB G$  are identical to those of the 
underlying pointed twisted curves. 

This can be seen a bit more explicitly, as follows. Let $P \to
\cC$ be the twisted $G$-cover associated to  $\cC \to \cB G$. As was
seen in Lemma \ref{Lem:P-representable}, the curve $P$ is projective,
and  deforming  $\cC \to \cB G$ is equivalent to deforming $P$ as a
$G$ space. Moreover, since  $P\to \cC$ is \'etale, this is
equivalent to deforming $\cC$ \cite{G-SGA1}\marg{Precise citation needed}.  

\subsubsection{\sc Obstructions:}
As in \cite{stablemaps}, Lemma 5.3.3, obstructions
of a pointed twisted curve $(\cC\to S, \Sigma_i^{\cC})$ are the same as
obstructions 
of the underlying $\cC \to S$, given by $\Ext^2(\Omega^1_\cC,
\cO_\cC)$. To show that this group vanishes, we follow 
\cite{Deligne-Mumford} closely. 

The local-to-global spectral sequence for $\Ext^*(\Omega^1_\cC,
\cO_\cC)$ 
involves, in degree 2, the three terms $$H^2(\cC,\cHom(\Omega^1_\cC,
\cO_\cC)), \quad  
H^1(\cC,\cExt^1(\Omega^1_\cC, \cO_\cC)),\quad\mbox{ and  }
H^0(\cC,\cExt^2(\Omega^1_\cC, 
\cO_\cC)).$$ We treat each term separately.
\begin{enumerate}
\item  By \cite{stablemaps}, Lemma 2.3.4, we have
$$H^2(\cC,\cHom(\Omega^1_\cC, \cO_\cC))=H^2(C,\pi_*\cHom(\Omega^1_\cC,
\cO_\cC))=0,$$ where $\pi:\cC\to C$ is the canonical morphism to the moduli
space. 
\item Similarly,
$\cExt^1(\Omega^1_\cC, \cO_\cC)$ is supported in dimension 0, hence the second
term vanishes. 
\item  We claim that the sheaf
$\cExt^2(\Omega^1_\cC, \cO_\cC))$ vanishes. This follows since locally in
$\cC$, the sheaf 
$\Omega^1_\cC$ has a 2-term locally free resolution. For instance, at a node
where the local picture for $\cC$ is $[U/\bmu_r]$ with $U=\Spec
k[z,w]/(zw)$, we have an  
$\bmu_r$-equivariant exact sequence on $U$:
$$ 0 \to \cO_U\stackrel{(z,w)}{\longrightarrow} \cO_U \oplus \cO_U
\stackrel{(dw, dz)}{\longrightarrow} \Omega^1_U \to 0,$$ with appropriate
$\bmu_r$-weights, giving a locally free resolution of $\Omega^1_\cC$.
\end{enumerate}
 Thus twisted curves are unobstructed.

\subsubsection{\sc Deformations:}
To calculate the dimension of  $\TC{g,n}{}{G}$ we  evaluate the
dimension 
of its tangent space.

We denote by $\cN_{\Sigma_i^{\cC}}$ the normal bundle of the $i$-th marking
$\Sigma_i^{\cC}\subset \cC$. 
As in \cite{stablemaps}, Lemma 5.3.2, we have an exact sequence
\begin{equation}\label{deformation-sequence}
A \to \Hom(\Omega^1_\cC, \cO_\cC) \to H^0(\cC, \oplus\cN_{\Sigma_i^{\cC}})
\to \Def 
\to 
\Ext^1(\Omega^1_\cC, \cO_\cC) \to 0,\end{equation}
where $\Def$ is the tangent space of the stack. 

{\em Infinitesimal automorphisms:} the space $A$ on the
left is the space of infinitesimal automorphisms of $(\cC,
\Sigma_i^{\cC})$, or, equivalently, of the twisted cover. We claim
that this space vanishes. This follows   
directly from the fact that $\TC{g,n}{}{G}$ is a Deligne-Mumford
stack, since the $\Isom$ schemes are unramified. 

We remark that this can also be 
computed on the level of twisted marked curves: 
 the sheaf 
of infinitesimal automorphisms is the subsheaf $\cA$ of
$\cHom(\Omega^1_{\cC}, \cO_{\cC})$ of  homomorphisms vanishing along
$\Sigma_i$. Considering its direct image in $C$, a local calculation
(similar to the one given below for the other terms)  reveals that it
is the same as the sheaf of infinitesimal automorphisms of $C$ fixing
the markings $\Sigma_i^C$, 
whose group of global sections
vanishes by the stability assumption. 

{\em The normal sheaf to a marking:}
at points where $\Sigma_i^{\cC}$ is untwisted we have that 
$H^0(\cC, \cN_{\Sigma_i^{\cC}})$ has dimension 1. At twisted markings
the 
normal space $ \cN_{\Sigma_i^{\cC}}$ has no nontrivial  sections: the
local picture of $\cC$ is 
$[U/\bmu_r]$, where $U=\Spec k[z]$ with
the standard action of $\bmu_r$, therefore $\bmu_r$ acts on a generator
$\partial/\partial z$ of $\cN_{\Sigma_i^{\cC}}$ via the nontrivial character
$\zeta_r 
\mapsto \zeta_r^{-1}$, and therefore the space of invariants is trivial.

{\em The extension groups:} the group $\Ext^i(\Omega^1_\cC, \cO_\cC)$
is dual to 
$H^{1-i}(\cC, \Omega^1_\cC\otimes \omega_\cC)$. By \cite{stablemaps},
Lemma 2.3.4 this is the same as $H^{1-i}(C, \pi_*(\Omega^1_\cC\otimes
\omega_\cC))$ 
where 
$\pi:\cC\to C$ is the natural map. Let us compare  $\pi_*(\Omega^1_\cC\otimes
\omega_\cC)$ with $\Omega^1_C\otimes \omega_C$. These are clearly
isomorphic away from the twisted markings and the twisted nodes. 

First consider a twisted
marking $\Sigma_i^{\cC}$ where the local picture of $\cC$ is the same as
 $[U/\bmu_r]$ with $U$ as above. The action of
$\bmu_r$ on $dz$ is 
via the standard character, therefore the invariant quadratic differentials are
generated by $z^r (dz/z)^2 = r^{-2} x (dx/x)^2 =r^{-2} (dx)^2/x$, where $x$ is
a 
parameter 
on $C$. That is, locally near such
a marking we have $\pi_*(\Omega^1_\cC\otimes
\omega_\cC) = \Omega^1_C\otimes \omega_C (\Sigma_i^C)$. 

Now consider a node on $\cC$. The local picture of $\cC$ is the same
as $[U/\bmu_r]$ with 
$U = \Spec k[z, w]/(zw)$, and the action can be described via $(z,w)
\mapsto (\zeta_r z, \zeta_r^a w)$ for some $a\in (\ZZ/l\ZZ)^{\times}$. The
sheaf $\omega_\cC$ has an invariant generator  $\nu_*(dz/z-dw/w)$, where $\nu$
is the normalization. The sheaf $\Omega^1_\cC$ has sections $f(z) dz + g(w) dw
+ \alpha z dw$. Invariant elements in $f(z) dz + g(w) dw$ are exactly
$\Omega^1_C/\operatorname{torsion}$, whereas $ z dw$ is invariant if and only
if $a=-1$, i.e. the node is balanced.

All in all we have $\chi(\pi_*(\Omega^1_\cC\otimes
\omega_\cC)) =\chi(\Omega^1_C\otimes \omega_C (\sum \Sigma_j^C) / \cT)$,
where the sum  $\sum \Sigma_j^C$ is taken over the {\em twisted} markings, and
the 
sheaf $\cT$ is the torsion subsheaf supported at {\em unbalanced}
nodes. 
It follows that $$\dim \Def = \chi(\pi_*(\Omega^1_\cC\otimes
\omega_\cC)) + H^0(\cC, \cN_{\Sigma_i^{\cC}}) =
H^0(C,\Omega^1_C\otimes \omega_C) + n - u.  $$ 
 The Proposition follows. \qed

The case $u=0$ corresponds to balanced twisted covers. We have:

\begin{corollary} The morphism $\TCbal{G}\to
\ocM_{g,n}$ is flat, proper and quasi-finite.
\end{corollary}

\proof  The morphism is proper and quasi-finite by Theorem
\ref{Th:stable-maps}, since $\bM=\bbS$, and therefore $\SM =
\ocM_{g,n}$. To check that it is flat it suffices to look at the map
of 
deformation spaces. But since the deformation spaces of source and
target are of the same dimension, and the map is quasifinite, it is
equidimensional. Since both are smooth, it follows from the local
criterion for flatness that the map is flat.\qed

\section{Twisted covers and admissible covers}\label{Sec:admissible}

In this section we compare our notion of twisted $G$-covers with the
notion of  admissible covers.

\subsection{Admissible covers} We recall the definition of an admissible cover
of nodal marked curves. Let $(C \to S, \Sigma_i)$ be an $n$-pointed nodal
curve of genus $g$, and let $d$ be a positive integer smaller than all the
residue 
characteristics of $S$.

\begin{definition} An admissible cover $p:D \to C$ of degree $d$ is a finite
morphism satisfying the following assumptions:
\begin{enumerate}
\item $D \to S$ is a nodal curve.
\item Every  node of $D$ maps to a node of $C$.
\item The restriction of $p:D \to C$ to $C\gen$  is
\'etale of constant degree $d$.
\item The local picture of $D \to C \to S$ at a point of $D$ over a
node of $C$ is the same as that of $D' \to C' \to S'$, with

$$\begin{array}{ccl}
D' & = &  \Spec A[\xi,\eta]/(\xi\eta-a)\\
\dar & &  \\
C' & = &  \Spec A[x,y]/(xy-a^e)\\
\dar & &  \\
S' & = &  \Spec A
\end{array}$$
for some integer $e\geq 0$, where $p^*x = \xi^e $ and $p^*y = \eta^e$.
\item For a geometric point over a marking of $C$, there is  an
integer $e\geq 0$ and an
analogous local picture
$$\begin{array}{ccc}
D' & = &  \Spec A[\xi]\\
\dar & &  \\
C' & = &  \Spec A[x]\\
\dar & &  \\
S' & = &  \Spec A
\end{array}$$
where $p^*x = \xi^e$, and $x$ is a local parameter for the marking.
\end{enumerate}
\end{definition}

We note that this generalizes the original definition of Harris and
Mumford in 
three ways:
\begin{itemize}
\item the genus of $C$ is arbitrary,
\item the ramifications over $\Sigma_i$ are not assumed to be simple, and
\item the curve $D$ is not required to have connected fibers.
\end{itemize}

Moduli of admissible covers in various degrees of generality were
discussed previously by Mochizuki \cite{Mochizuki}, 
Wewers \cite{Wewers}, and \cite{A-O:Hurwitz}. See also
Bernstein \cite{Bernstein}, Bouw--Wewers \cite{Bouw-Wewers}.

 For the rest of the
section we fix $\bbS =  \Spec \ZZ[1/d!]$.
Admissible covers of degree $d$ of {\em stable} $n$-pointed curves of
genus $g$ 
form a proper Deligne--Mumford stack $\Adm_{g,n,d}\to \Spec \ZZ[1/d!]$
admitting a projective coarse moduli space, see \cite{Mochizuki}. See
also  \cite{A-O:Hurwitz}. 


\subsection{Twisted covers and the normalization of the
Harris--Mumford stack} 
Since the appearance of \cite{Harris-Mumford}, there has been some
dissatisfaction with the stack of admissible covers, for two
reasons. First, the original definition involves a description of
{\em families} of admissible  covers, as the moduli problem was not
determined by the geometric objects it parametrized - a resolution of
this issue using logarithmic structures is given in Mochizuki's work
\cite{Mochizuki}; our approach below uses twisted curves instead. Second,   
it follows from the description of the deformation spaces of
admissible covers 
in \cite{Harris-Mumford}\marg{Precise citation needed} that $\Adm_{g,n,d}$ is
in general not 
normal, 
but its normalization is always smooth. Below we exhibit this
normalization as 
a stack of twisted covers.

We build on the usual equivalence of categories 
$$\left\{ \mbox{\parbox{2.8in}{\begin{center}finite \'etale
covers $D\to S$
of degree $d$\end{center}}}\right\} \leftrightarrow \left\{
\mbox{\parbox{2.1in}{\begin{center}principal $\cS_d$-bundles $P\to
S$\end{center}}}\right\},$$ 
 where $\cS_d$ is the symmetric group on $d$
letters. We briefly recall that to a principal $\cS_d$-bundle $P \to S$ one
associates the finite \'etale cover $D\to S$ where $D = P/\cS_{d-1}$;
and given 
a finite \'etale cover $D \to S$ of degree $d$, we consider the
complement $P$ of all the 
diagonals in the $d$-th fibered power $D^d_S$, which is easily seen to be a
principal $\cS_d$-bundle.

Consider a {\em balanced} twisted $\cS_d$ cover  $P \to \cC$. The
schematic quotient $D = P/\cS_{d-1}$ is not necessarily \'etale over
$C$. Instead we have the following:


\begin{lemma}\label{Lem:quotient-is-admissible} The morphism $D\to C$ 
is an admissible cover of degree $d$. 
\end{lemma}
\proof Indeed $D = P/\cS_{d-1}$ is a nodal curve, being the quotient
of a nodal curve by a group acting along the fibers. Since every node
of $D$ is the image of a node on $P$, we have that its image in $C$ is
  again a node. Since $P\gen \to C\gen$ is a principal $\cS_d$-bundle, its
quotient $D\gen$ by $\cS_{d-1}$ is finite and \'etale of degree
$d$. Finally, given a geometric point of $D$ over a node of $C$, we
can choose a point $z'$ in the preimage in $P$, which is nodal. The
local picture is the same as $\Spec A[u,v]/(uv-b))$. The
stabilizer in $\cS_{d}$ of $p$ is a cyclic group $C_r$ which acts via $(u,v
)\mapsto (\zeta_ru, \zeta_r^{-1}v)$, in such a way that $x = u^r$ and
$y=v^r$.  If $C_{r'} = C_r\cap
S_{d-1}$, write $\xi = u^{r'}$ and $\eta = v^{r'}$, and the local
picture of $D$ is $\xi\eta = b^{r'}$. Setting $e = 
r/r'$, 
it follows that  $x
= \xi^e$ and  $y = \eta^e$, as required. A similar (and simpler)
argument gives the structure along a marking.\qed

 The functor, which associates to a twisted $\cS_n$-cover $P \to \cC$
the admissible cover $P/\cS_{n-1}\to C$, is a morphism of stacks
$\phi:\TCbal{\cS_d} \to \Adm_{g,n,d}$.

\begin{proposition}
 The morphism $\phi$ exhibits  $\TCbal{\cS_d}$ as the normalization of
 $\Adm_{g,n,d}$. 
\end{proposition}

\proof It suffices to show that $\phi$ is finite and surjective, and
there is 
an open dense substack $\Adm_{g,n,d}^0$, whose inverse image is dense
in $\TCbal{\cS_d}$, over which $\phi$ is an isomorphism. 

The morphism $\phi$ is finite since 
\begin{enumerate}
\item $\phi$ is representable: for this it suffices to show (see, e.g.,
\cite{stablemaps}, Lemma 4.4.3) that any automorphism of a twisted
$\cS_d$-bundle 
$P\to \cC$
which acts trivially on the associated admissible cover $D \to C$ is the
identity. Indeed, such an isomorphism acts trivially on $D\gen \to C\gen$, and
by the equivalence of categories cited above it acts trivially on $P\gen$ as
well. By \cite{stablemaps}, Theorem 4.4.1 (and also Remark 4.4.4)  an
automorphism of a twisted stale map 
$\cC \to \cM$ is
determined by its action on the generic object $C\gen \to \cM$, and the claim
follows.

\item $\phi$ is proper, since $\TCbal{\cS_d}$ is proper.
\item $\phi$ is quasifinite, since, by Theorem \ref{Th:stable-maps},
the stack  $\TCbal{\cS_d}$ is
quasifinite over $\ocM_{g,n}$.
\end{enumerate}

Consider the dense open  substack $\Adm_{g,n,d}^0$ of admissible
covers of {\em smooth} 
curves $C$. Its inverse image in $\TCbal{\cS_d}$ is the dense open
substack of balanced twisted covers over smooth twisted curves. We
claim that over $\Adm_{g,n,d}^0$  
the morphism $\phi$ is 
an 
isomorphism. It is easy
to see that $\Adm_{g,n,d}^0$ is smooth, therefore we need only construct an
inverse functor of $\phi$ for admissible covers over {\em reduced} base
schemes. Consider an 
admissible cover $D \to C$ with $C\to S$ smooth and 
$S$ reduced. Then $D$
is a smooth curve branched only over the marked sections  of $C$. The
restriction $D\gen\to C\gen$ is a finite \'etale cover corresponding to a
principal $\cS_d$-bundle $P\gen \to C\gen$. The tameness assumption and
Abhyankar's lemma (see \cite{G-SGA1}, exp. XIII section 5) imply that
the normalization $P$ of $C$ in 
the structure sheaf of $P\gen$ is again a smooth curve, the quotient $\cC
=[P/\cS_d]$ is a smooth 
twisted pointed curve over $S$, and $P \to \cC$ is a twisted $\cS_d$-bundle. 
This
provides an inverse of $\phi$ restricted over $\Adm_{g,n,d}^0$. 

 Next we show that $\phi$ is surjective. We give two arguments for this fact,
since we feel they are instructive in different ways.

\subsubsection{Surjectivity I} For the first argument, consider an admissible cover $D_0 \to C_0$
over an algebraically closed 
field. In \cite{Harris-Mumford} it is shown that its deformation space is
reduced, and the locus of smooth admissible covers in it is nonempty. Thus
$\Adm_{g,n,d}^0$ is everywhere dense, and by properness $\phi$ surjects on a
dense closed substack. But since $\Adm_{g,n,d}$ is reduced, $\phi$ is
surjective.

\subsubsection{Surjectivity II} The second argument is longer, but more elementary, as it does not
use deformation 
theory. Let $D \to C$ be an admissible cover defined over an
algebraically closed field; we will produce a twisted cover of $C$
according to 
the 
following procedure. We think of $\cB\cS_d$ as the stack of \'etale
covers of 
degree 
$d$; the restriction of $D$ to $C\gen$ gives a generic object $C\gen \to
\cB\cS_d$. 
Let us produce charts for this object in the sense of
\cite{stablemaps}, Section 3.2.

Let $p \in C$ be a marked point, and let $m$ be the least common
multiple of all 
the ramification indices of points of $D$ over $C$. Let $U \to C$ be a
morphism from a 
smooth, but not necessarily complete, curve $U$ such that

\begin{enumerate}

\item the image of $U$ does not contain any special point of $C$ except $p$,

\item there is precisely one point $q\in U$ over $p$, and

\item there is an action of a cyclic group $\Gamma$ of order $m$ on $U$, having
$q$ 
as a fixed point and leaving the morphism $U \to C$ invariant, which is free
outside of $q$ and such that the induced morphism $U/G \to C$ is \'etale.

\end{enumerate}

The normalization $\widetilde D$ of the pullback of $D$ to U is \'etale over
$U$, and the action of $G$ on $U$ lifts to an action on $\widetilde D$; this
gives 
a chart around the point $p$. It is easily checked that the quotient
$\widetilde 
D/\Gamma$ is the pullback of $D$ to $U/\Gamma$.

If $p\in C$ is a node, the procedure is similar. Let $m$ be the least common
multiple of all the ramification indices of points of $D$ over $C$.
Let $U \to C$ a morphism from a nodal, but not necessarily complete, curve
$U$ such that

\begin{enumerate}

\item the image of $U$ does not contain any special point of $C$ except $p$,

\item there is precisely one point $q\in U$ over $p$,

\item $U$ has two irreducible components $U_1$ and $U_2$, which are smooth and
intersect only at $q$,

\item there is a balanced action of a cyclic group $\Gamma$ of order $m$ on
$U$, 
having $q$ as a fixed point and leaving the morphism $U \to C$ invariant, which
is 
free outside of $q$ and such that the induced morphism $U/G \to C$ is \`etale.
\end{enumerate}

Recall that the action of $\Gamma$ on $U$ is {\em balanced} when the two
characters 
of $\Gamma$ describing the action of $\Gamma$ on the tangent spaces to $U_1$ and
$U_2$ are opposite.

Let $\widetilde D_1$ and $\widetilde D_2$ be the normalizations of the
pullbacks of 
$D$ to $U_1$ and $U_2$ respectively; then $\widetilde D_i$ is smooth over
$U_i$. 
The action of $\Gamma$ on $U_i$ lifts to an action on $\widetilde D_i$. To
obtain a 
chart, we choose a way of identifying the fiber of $\widetilde D_1$ over $q$
with 
the fiber of $\widetilde D_2$ over $q$; this gives a $\Gamma$-equivariant
\`etale 
cover $\widetilde D \to U$, and the quotient $\widetilde D/\Gamma$ is precisely
the 
pullback of $D$ to $U/\Gamma$. This gives the desired chart.\qed
\marg{In the final version there should be a subsection about the
singularities of $\phi:\TCbal{\cS_d} \to \Adm_{g,n,d}$.*}

\subsection{Admissible $G$-covers and twisted covers}\label{Sec:admissible-G-covers}

\begin{definition}
Let $C$ be a nodal curve over a scheme $S$. An admissible
$G$-cover $\phi:P \to C$ is an admissible cover, with an action of $G_S$ on
$P$ 
leaving $\phi$ invariant, satisfying the following two conditions.

\begin{enumerate}

\item The restriction $P\gen \to C\gen$ is a principal $G_S$-bundle.

Let $p$  be a geometric point of $P$. Notice that this first
condition 
insures that the stabilizer $G_p$ of $p$ is a cyclic group. Then we also assume

\item for each geometric  nodal point $p$ of  $P$, with image point
$s$ of the base
 $S$, the action of the stabilizer
$G_p$ of $p$ 
on the fiber $P_s$ over $s$ of $P \to S$  is balanced.

\end{enumerate}
\end{definition}

Admissible $G$-covers form a category
$\Adm_{g,n}(G)$ with arrows given by fiber diagrams.

In contrast with the case of plain admissible cover, we have the following
result. 

\begin{theorem}
There is a base-preserving equivalence of categories between $\TCbal{G}$ and
$\Adm_{g,n}(G)$.
\end{theorem}

So in particular $\Adm_{g,n}(G)$ is a Deligne--Mumford stack, isomorphic to
$\TCbal{G}$.  \marg{In final version this should be proven directly in
appendix.} 
See also \cite{Wewers}, \cite{Bouw-Wewers}.  

\begin{proof}
Let $\cC$ be a twisted curve over a scheme $S$, with moduli space $C$. If $P
\to 
\cC$ is a balanced twisted $G$-cover, we have seen in
\ref{Sec:twisted-G-covers} that $P$ is an algebraic space.
 We claim that the
composition of 
$P\to \cC \to C$, with $\cC \to C$ the morphism to the moduli space, is an
admissible $G$-cover. The fact that $P\gen \to C\gen$ is a principal
bundle follows from the definition, since $C\gen = \cC\gen$ and $P \to
\cC$ is a principal bundle. The fact that $P \to C$ is an admissible
cover, as well as condition (2) above, are identical to the argument
of Lemma \ref{Lem:quotient-is-admissible}.

Conversely, given an admissible $G$-cover $P\to C$ over a scheme $S$, consider
the stack quotient 
$\cC = [P/G_S]$. Now $P\to \cC$ is a principal $G$-bundle, the morphism $\cC
\to \cB G$ is representable since $P$ is, $\cC$ is nodal being the quotient of
a nodal curve, $C$ is the moduli space of $\cC$, and  the action of $G$ on $P$
is balanced, showing that $\cC$ is balanced.

It is easy to see that these correspondences are functorial, and that they are
inverse to each other in the usual sense.
\end{proof}

\section{Rigidification and Teichm\"uller structures}

In this section we define the notion of {\em rigidification} of a
stack, and use it to define {\em twisted Teichm\"uller $G$-strauctures}.

\subsection{Rigidification}\hfill

\subsubsection{The example of line bundles} The classic instance of
rigidification of a stack occurs when constructing the Picard scheme
of an irreducible projective variety. 

Consider the stack $\cB
\GG_m$, the classifying stack of the multiplicative group. 
As it stands, its
objects are $\GG_m$-bundles, but it is well known to be equivalent to
the 
stack of line bundles: objects over a scheme $X$ are line bundles over
$X$, and a morphisms between $L_1 \to X_1$ and $L_2 \to X_2$  are
fiber squares 
$$ \begin{array}{ccc} L_1 & \to & L_2 \\
\dar & \square  & \dar \\ 
X_1 & \to & X_2.
\end{array}
$$

Of course every line bundle over $X$ has the group of invertible
functions $\cO_X(X)^\times$ in its automorphism group, and in
particular, if $X$ is a projective variety over a field $k$, the automorphism
group contains $k^\times$. For such $X$  we can define the
``Picard stack'', whose objects over a $k$-scheme $T$ are line bundles
on $T\times X$. Clearly this stack is not representable, since  every 
geometric object has automorphisms. 
How do we obtain the Picard {\em scheme}
of, say,
an irreducible projective  variety $X$ out of $\cB \GG_m$?

The point is, that one can take the category of line bundles over $X$
and {\em rigidify} it by removing the multiplicative group from the
automorphisms of any object. The traditional
procedure (see \cite{G-FGA}) is to take, as a first  
approximation, the category whose objects are line bundles over
$T\times X$ and whose arrows are ``isomorphisms up to twisting'' $L_1
\to L_2 
\otimes M$, where $M$ is a line bundle coming from $T$. This works
\'etale-locally, and one needs to sheafify the
resulting category in the \'etale topology. The result is precisely
the Picard scheme. In 
effect, this procedure produces a stack whose geometric objects are
still line bundles over $X$, but where the automorphisms given by the
multiplicative group are ``removed from the picture''.

\subsubsection{Level structures and rigidified covers}

Suppose one is interested in studying full level-$m$ structures on
smooth curves of genus $g>1$. A structure of full level-$m$ on a smooth
curve $C$ of genus $g$ over an 
algebraically closed field is a basis for the
$\Zm$-module $\opH^1(X,\Zm)$. This basis corresponds to an element of
$H^1\bigl(X,\level\bigr)$, with the property that the associated
$\level$-principal bundle $E \to C$ is connected. One would be temped to
identify this level structure with the principal bundle $E \to C$,
thought of as an element of 
$\TC{g,0}{}{\level}$, but this would be an error, because the group
$\level$ acts on the bundle $ E \to C$, while a level
structure should have no nontrivial automorphisms fixing $C$. 

Considering twisted $G$-covers in general, the center $Z(G)\subset G$ 
always acts on any twisted 
$G$-cover; we wish to {\em rigidify} the twisted covers by making this action
trivial. We now describe a procedure for this in some generality.


\subsubsection{Rigidification in terms of a presentation}
The idea of rigidification can be seen explicitly with a
presentation. 

Let 
$H$ be a flat finitely presented separated group scheme over  
a base  scheme $\bbS$, and
${\cX}$ an algebraic  stack over $\bbS$.  
Take a smooth map of finite presentation $U \to {\mathcal X}$, and set $R
= U \times_{\mathcal X}U$, so that $R \double U$ is a smooth presentation
for ${\mathcal X}$. 

Assume that there
is an action of $H_{U\times_SU}$ on $R$, leaving the two projections $R
\to U$ invariant.  Then there exists a quotient smooth
groupoid $R/H \double U$; this is a smooth presentation of a  stack
${\cX^H}$, which is the rigidification of $\cX$, removing $H$ from the
stabilizers.

We now describe a natural situation where such a picture occurs.

\subsubsection{The rigidification setup}

Again, let $H$ be a flat finitely presented separated group scheme over 
a base  scheme $\bbS$,
${\cX}$ an algebraic  stack over $\bbS$. Assume that for each object
$\xi \in 
{\cX}(S)$ there is an embedding $$\iota_\xi\colon H(S)\into \Aut_S
(\xi),$$ which is compatible with pullback, in the following sense:
given two objects $\xi\in {\cX}(S)$ and $\eta\in{\cX}(T)$, 
and an arrow $\phi\colon \xi \to \eta$ in ${\cX}$ over a 
morphism of schemes $f\colon S \to T$, the natural pullback homomorphisms
$$\phi^*\colon \Aut_T (\eta)\to \Aut_S (\xi)$$ and $$f^* \colon H(T) \to H(S)$$
commute with the embeddings, that is, $\iota_\xi f^* = \phi^*\iota_\eta$.

This condition can also be expressed as follows. Let $\phi\colon
\xi \to \eta$ be an arrow in ${\cX}$ over a morphism of schemes
$f\colon S \to
T$, and
$g \in H(T)$. Then
the diagram
$$\begin{array}{ccc}
	 \xi            & \mapright \phi   & \eta      \\
         \mapdown{f^*g} &                  &\mapdown g \\
         \xi            & \mapright \phi &\eta       
\end{array}
$$
commutes. In particular, by taking $\xi = \eta$ and $\phi$ to be in
$\Aut_S(\xi)$, we see that $H(S)$ must be in the center of
$\Aut_S(\xi)$. In particular, $H$ must be commutative. (One might
consider a more general situation,  where the element $g$ on the right
differs from the element appearing on the left, but this is not
crucial for our purposes.)

\begin{theorem}
Let $\cX\to \bbS$ be an algebraic stack, $H\to \bbS$ a
flat finitely presented group 
scheme over $\bbS$, and assume that for every object $\xi\in
\cX(T)$ there is an embedding $H_T\subset \Aut_T(\xi)$ compatible with
pullbacks. Then there is a smooth surjective finitely presented morphism
of
algebraic stacks
$\cX
\to \cX^H$ 
satisfying the following properties:
\begin{enumerate}
\item For any object $\xi\in \cX(T)$ with image $\eta\in \cX^H(T)$, we
have that $H(T)$ lies in the kernel of $\Aut_T( \xi) \to \Aut_T(\eta)$.
\item The morphism $\cX \to \cX^H$ is universal for  morphisms of stacks
$\cX \to \cY$
satisfying (1) above.
\item If $T$ is the spectrum of an algebraically closed field, then in
(1) above, $\Aut_T(\eta) = \Aut_T( \xi)/H(T)$.
\item A moduli space for $\cX$ is also a moduli space for $\cX^H$.
\end{enumerate}

Furthermore, if $\cX$ is a Deligne--Mumford stack, then $\cX^H$ is
also a Deligne--Mumford stack and the morphism $\cX \to \cX^H$ is
\'etale.

\end{theorem}

People familiar with the theory of $n$-stacks might recognize this
rigidification as the stack associated to the quotient $[\cX/\cB H]$
of $\cX$ by the action of the group-stack $\cB H$. Out proof below
takes a slightly more concrete view.

\subsubsection{The action on $\underhom$-sheaves}

Given an object $\xi$ of $\cX$ over a scheme $S$, the embeddings $\iota_\xi$
define a categorically injective morphism of $S$-group schemes of the
pullback
$H_S$ of
$H$ to
$S$ to the group scheme $\underaut_S(\xi)$ of automorphisms of $\xi$.

With these hypotheses, if $\xi$ and $\eta$ are objects of $\cX$ over
two schemes $S$ and $T$ respectively, and $f \colon S \to T$ is a morphism
of schemes, there is an action of $H(T)$ on the set $\Hom_f(\xi,\eta)$ of
arrows in ${\cX}$ lying over $f$, defined by setting $g \cdot \phi =
g\circ \phi = \phi f^*g$ for each $g \in H$ and each
$\phi\in\Hom_{\cX}(\xi,\eta)$. If $f\colon S \to T$ and $f' \colon T
\to U$ are morphisms of schemes, and $\xi$, $\eta$ and $\zeta$ are object
of $\cX$ over $S$, $T$ and $U$ respectively, there is a composition map
$$
\Hom_f(\xi,\eta) \times \Hom_{f'}(\eta,\zeta) \longrightarrow
\Hom_{f'f}(\xi,\zeta);
$$
it is easy to see that map passes to the quotient, yielding a map
$$
\Hom_f(\xi,\eta)/H(T) \times \Hom_{f'}(\eta,\zeta)/H(U) \longrightarrow
\Hom_{f'f}(\xi,\zeta)/H(U).
$$

This action of $H(T)$ on $\Hom_f(\xi, \eta)$ induces a right
action of the group scheme $H_T$ on the sheaf $\underhom_f(\xi, \eta)$,
which sends each scheme $T'$ over $T$ to the set of arrows in $\cX$ from
the pullback of $\xi$ to $S \times_T T$ to the pullback of $\eta$ to $T'$
lying over the projection $S \times_T T'  \to T'$.

\subsubsection{Rigidification in categorical terms}
Consider the quotient sheaf 
$$\underhom_f^H(\xi, \eta) = \underhom_f(\xi,
\eta)/H_T \colon ({\rm Sch}/T)^{\rm opp} \to ({\rm Sets}),$$ that is, the
quotient sheaf associated to the presheaf which sends each $T' \to T$ to
the set $\Hom_f(\xi, \eta)(T')/H(T')$. We define $\Hom_f^H(\xi,
\eta) = \underhom_f(\xi, \eta)(T)$ to be the set of global sections of
this sheaf. The composition map above induces a morphism of sheaves of
sets on $U$
$$
f'_*\underhom_f^H(\xi,\eta) \times
\underhom_{f'}^H(\eta,\zeta) \longrightarrow
\underhom_{f'f}(\xi,\zeta)
$$
and hence a map
$$
\Hom_f^H(\xi,\eta) \times \Hom_{f'}^H(\eta,\zeta) \longrightarrow
\Hom_{f'f}^H(\xi,\zeta).
$$

We define a category ${\cX}^H_{\pre}$, in which the objects are
objects of $\cX$, and an arrow from $\xi$ to $\eta$ consist of an
element of $\Hom_f^H(\xi,\eta)$ for some
$f \colon S \to T$, where $S$ and $T$ are the schemes
underlying $\xi$ and 
$\eta$ respectively, and the composition is defined by the map above. This
${\cX}^H_{\pre}$ is a prestack(\cite{L-MB}, D\'efinition~3.1), but not a
stack, in general; we
define the category ${\cX}^H$ to be the stack associated to
${\cX}^H_{\pre}$, as in \cite{L-MB}, Lemme~3.2. We note that the
process of taking the stack associated to a prestack has the property
that ${\cX}^H(T) = {\cX}^H_{\pre}(T)$ whenever $T$ is the spectrum of
an algebraically closed field.

There are obvious functors from $\cX$ to ${\cX}^H$ and from ${\cX}^H$
to the category of schemes over $\bbS$. It follows  from the 
construction that if $T$ is the spectrum of an algebraically closed field,
and $\xi$ and $\eta$ are objects of ${\cX}$, then the set of
isomorphisms of $\xi$ and $\eta$ in ${\cX}^H$ is the set of
isomorphisms of $\xi$ and $\eta$ in ${\cX}(T)$ divided by the natural
action of $H(T)$.

It is easily checked that ${\cX}^H$ is a stack fibered in groupoids over
$\bbS$, and that it has properties~(1), (2) and (3) of the theorem.
Property~(4) follows immediately from property~(3).

We claim that ${\cX}^H$ is in fact an algebraic stack. First of all let us
show that the diagonal of ${\cX}^H$ is representable, finitely presented and
separated. Let $X$ and $Y$ be schemes,
$X \to {\cX}^H$ and $Y \to {\cX}^H$ two morphisms corresponding to objects
$\xi 
\in {\cX}^H(X)$ and $\eta \in {\cX}^H(Y)$. We need to show that the fiber
product
$X
\times_{{\cX}^H}Y$ is  representable, separated and of finite
presentation over
$X
\times_\cS Y$. This is a local question in the flat topology over $X$ and
$Y$, so we may suppose that
$\xi$ and $\eta$ are objects of $\cX(X)$ and $\cX(Y)$ respectively. Then $X
\times_{{\cX}^H}Y$ corresponds to the
functor
$\underisom_{X
\times_{\bbS} Y}^{{\cX}^H}({\rm pr}_X^* \xi, {\rm pr}_Y^*\eta)$, which is by
definition equal to the quotient of the algebraic space $X
\times_{\cX} Y = \underisom_{X 
\times_{\bbS}Y}^{{\cX}}({\rm pr}_X^* \xi, {\rm pr}_Y^*\eta)$ by the
action of $H_U$.  It is easy to see that
this action is free, in the sense that the morphism
$$
\underisom_{X \times_{\bbS}Y}^{{\cX}}({\rm pr}_X^* \xi, {\rm pr}_Y^*\eta)
\times_T H_T \longrightarrow 
\underisom_{X \times_{\bbS}Y}^{{\cX}}({\rm pr}_X^* \xi, {\rm pr}_Y^*\eta)
\times_T \underisom_{X \times_{\bbS}Y}^{{\cX}}({\rm pr}_X^* \xi, {\rm
pr}_Y^*\eta)
$$
which is the action on one component and the projection on the other is
categorically injective. By a result of M. Artin (\cite{L-MB},
Corollarie~10.4.1), the quotient 
$X \times_{{\cX}^H} Y$ is an algebraic space over $X \times_\cS
Y$. It is easily checked that it is separated and of finite presentation.

Next we'll prove that the morphism $\cX \to \cX^H$ is smooth, surjective
and of finite presentation. If $X \to \cX$ is a smooth surjective map
locally of finite presentation from a scheme $X$, we claim that the
composition $X \to 
\cX  \to \cX^H$ is smooth, surjective and locally of finite
presentation, which implies that $\cX^H$ is algebraic.

To check this, it is enough to show that given a morphism $T\to
\cX$, where
$T$  is a scheme, the fiber product $\cX \times_{\cX^H} T$
smooth. surjective and of
finite presentation on $T$. By passing to a flat cover of
$X$, we may assume that this morphism $T \to \cX^H$ factors through $\cX$;
in other words, it is enough to show that the projection $\cX
\times_{\cX^H}
\cX \to \cX$ is smooth, surjective and of finite presentation. An object
of
$\cX \times_{\cX^H} \cX$ is given by a scheme $T$ over $\bbS$, two objects
$\xi$ and $\eta$
of $\cX$ over $T$, and an isomorphism $\alpha$ of $\xi$ with $\eta$ in
$\cX^H(T)$. Consider the principal $H$-bundle $\underisom_T^{\cX}(\xi,
\eta)
\to \underisom_T^{\cX^H}(\xi, \eta)$; by pulling it back to $T$ via the
morphism
$T\to \underisom_T^{\cX^H}(\xi, \eta)$ given by $\alpha$ we get a principal
$H$-bundle $P \to T$. There is an obvious functor $\cX \times_{\cX^H}
\cX \to \cX \times_\bbS \cB_\bbS H$ sending the object $( \xi, \eta, \alpha)$
to to $(\xi, P)$. I claim that this is an isomorphism. Let us define the
inverse functor $\cX \times_\bbS \cB_\bbS H \to \cX \times_{\cX^H}
\cX \to \cX$. Take an object $(\xi, P)$ of $\cX \times_\bbS \cB_\bbS H$ over
a scheme $T$, and consider the pullback $\xi_P$ of $\xi$ to $P$; the
embedding $H_P \into \underaut_P(\xi)$ defines an action of $H_T$ on
$\xi_P$, giving descent data to descend $\xi_P$ to another object $\eta$ of
$\cX$ over $T$. By definition, this $\eta$ comes equipped with a canonical
isomorphism $\alpha\colon \xi\simeq \eta$ in $\cX^H(T)$. We define the
image of the object $(\xi, P)$ to be the object $(\xi, \eta , \alpha)$ of
$\cX \times_{\cX^H} \cX$. We leave it to the reader to define the action of
this fuctor on arrows, and to check that it gives an inverse to the functor
above.

So the projection $\cX \times_\bbS\cX \to \cX$ is isomorphic to the
projection $\cX \times_\bbS \cB_\bbS H \to \cX$. Obviously $ \cB_\bbS H$
is surjective and of finite presentation ovewr $\bbS$; we only need to
check that it is smooth. This is obvious when $H$ is smooth over
$\bbS$; in general the morphism $S \to  \cB_\bbS H$ given by the trivial
torsor is flat and surjective, but not smooth. The result follows from the
following lemma.

\begin{lemma}
Let $\cX$ be an algebraic stack flat of finite presentation over a scheme
$S$, and assume that there exists a flat surjective morphism $U \to \cX$,
where $U$ is smooth over $S$. Then $\cX$ is smooth over $S$.
\end{lemma}

\begin{proof}
The statement is local in the smooth topology on $\cX$, so we may assume
that $\cX$ is a scheme; in this case the result is standard.
\end{proof}

To conclude the proof of the theorem, we assume that $\cX$ is
Deligne--Mumford, and show that the morphism $\cX \to \cX^H$ is \'etale.
As we saw above, it is enough to check that the projection $\cX
\times_{\cX^H} \cX \to \cX$ is \'etale. Take a morphism $T \to \cX$, where
$T$ is a scheme; we have seen that the fiber product $\cX \times_{\cX^H}
T$ is isomorphic to $T \times_{\bbS}\cB_\bbS H = \cB_T H_T$, so we only
need to prove that $\cB_T H_T$ is \'etale over $T$, or, equivalently, that
$H_T$ is \'etale over $T$. But this is clear, because if $\xi\in \cX(T)$
is the object corresponding to the given morphism $T \to \cX$, then there
is an embedding $H_T \into \underhom_T(\xi)$, and $\underhom_T(\xi)$ is
unramified over $T$, because $\cX$ is Deligne--Mumford.

Furthermore, since any morphism from $\cX$ to an algebraic space
factors uniquely though ${\cX}^H$, we see that $\cX$ and ${\cX}^H$
share the same moduli space. \qed

\begin{definition} Let $\cX $ and $  H$ be as in the theorem. Then
we call $\cX \to \cX^H$ the {\em rigidification of $\cX$ along
$H$}. Objects of $\cX^H$ are called 
{\em $H$-rigidified objects of $\cX$.} 
\end{definition}


\subsection{Teichm\"uller  structures}

In this section we assume $n=0$

\begin{definition} Denote $\TC{g}{\rig}{G} =  (\TC{g}{}{G})^{Z(G)}$,
the stack of $Z(G)$-rigidified twisted $G$-covers. Let $\TCtei{G}\subset
\TC{g}{\rig}{G}$ be the open-and-closed substack whose geometric
objects correspond to  {\em connected}, balanced, rigidified twisted
$G$-covers. We call this the stack of {\em twisted Teichm\"uller
$G$-structures}. We denote by  $\TCtei{G}^0$ the open substack  of smooth
Teichm\"uller 
$G$-structures, namely twisted Teichm\"uller structures over smooth curves. 
\end{definition}

\begin{lemma}
The
morphism $\TCtei{G}^0\to \cM_g$ is finite and \'etale. 
\end{lemma}
\proof This morphism is \'etale since the deformation space of a  smooth
Teichm\"uller 
structure coincides with that of the underlying curve. We need to show that
this morphism is representable, which, by a well known result (see
\cite{stablemaps}, Lemma 
4.4.3) means that  
the induced maps of automorphism groups of geometric objects is injective. But
the automorphism group $\Aut_C P$ of a connected 
principal bundle $P \to C$ over an algebraically closed field  fixing $C$
is the center of the structure group, therefore, when we rigidify the
automorphism group becomes trivial. \qed

In \cite{Pikaart-DeJong}, a stack of ``Techm\"uller $G$-level
structures over smooth curves''
${}_G\cM_g$ is defined for any finite group $G$. This generalizes the
treatment of Deligne and Mumford in \cite{Deligne-Mumford}, Section
5. For the benefit of 
the reader familiar with their construction, we compare it with our
stack. We note that this result is not necessary for understanding the
rest of this paper.

\begin{proposition}
Assume $G$ is a constant finite group. The stack
 $\TCtei{G}^0$ is isomorphic to the stack ${}_G\cM_g$ of Teichm\"uller level
 structures.
\end{proposition}
\proof Given an object of $\TCtei{G}^0(S)$, then \'etale locally on the base
$S$ we have a principal $G$-bundle $P \to C$ as well as a section $s:S \to C$.
By definition this yields a surjective homomorphism $\pi_1(C/S,s) \to G$. It is
straightforward to verify that this 
only depends on the original object and yields a global section of the sheaf
$\cHom^{\operatorname{ext}}(\pi_1(C/S), S)$, namely an object of
${}_G\cM_g(S)$. It is also easy to check that this is functorial.

 Both
stacks $\TCtei{G}^0$ and  ${}_G\cM_g$ are finite and \'etale over
$\cM_g$. In order to check that this functor gives an isomorphism it suffices
to check that it is bijective on geometric points. Now, if $C$ is a curve over
an algebraically closed field, then there is a one-to-one correspondence
between isomorphism classes of connected $G$-covers and surjective
homomorphisms $\pi_1(C,s) \to G$ up to conjugacy. \qed

Pikaart and De Jong (again generalizing Deligne and Mumford) proceed to define
a proper stack $\overline{{}_G\cM_g}$ by 
normalizing $\overline{\cM_g}$ in  ${}_G\cM_g$. There is still a functor
$\TCtei{G}\to \overline{{}_G\cM_g}$, which is in general not an
isomorphism: for instance, $\overline{{}_G\cM_g}$ is in general singular, and
the morphism $\overline{{}_G\cM_g}\to \overline{\cM_g}$ is always
representable. In contrast, $\TCtei{G}$ is always nonsingular,
and in general $\TCtei{G}\to \overline{\cM_g}$ is not
representable\marg{In the final version some examples of these
singularities are due}.

\section{Abelian twisted level structures}

{\bf Convention.} Fix a positive integer $m$, and set $\bbS = \Spec \ZZ[1/m]$.
 Throughout this section $G = (\ZZ/m\ZZ)^{2g}$.

In this section we study in detail the case $G = (\ZZ/m\ZZ)^{2g}$. By
elementary covering theory, a
smooth Teichm\"uller $G$-structure on a curve $C$ 
consists of a basis for $H^1(C,\ZZ/m\ZZ)$, which is what often one
calls a level-$m$ curve. Below we extend this
description to stable curves: a twisted Teichm\"uller $G$-structure
over a base $S$ is
equivalent to a {\em  twisted level-$m$}  curve
\ref{Th:twisted-level-m}, which is a {\em pre-level-$m$} twisted curve
 $h:\cC \to S$ along 
with a basis of  $\bR^1h_*\ZZ/m\ZZ$ which is a local system.

\subsection{Pre-level-$m$ curves}

\begin{definition}
A balanced twisted nodal curve $\cC$ whose moduli space $C$ is stable  is said
to be  a 
{\em pre-level-$m$ curve} if for each geometric fiber,
the stabilizer at each separating node  is trivial and
the stabilizer at a 
non-separating node is cyclic of order $m$.
\end{definition}

\begin{proposition}\label{Prop:pre-level}
The underlying twisted curve of an
object of   $\TCtei{\level}$ is a pre-level-$m$ curve.
\end{proposition}

\proof By definition, it is enough to consider a geometric object, and such an
object is the rigidification of a connected twisted $G$-cover $P \to \cC$. The
claim is 
obvious for a smooth curve, therefore we may assume $\cC$ is nodal.
The stack of unpointed balanced twisted $G$-covers is quasifinite over
$\ocM_g$ and of pure dimension $3g-3$, therefore we can smooth the nodes
independently, in particular, for each node $z$ of $\cC$  the  twisted
$G$-cover $P \to \cC$  deforms to a curve $P_\eta \to \cC_\eta$ with exactly
one node $z_\eta$ having $z$ in its closure. Since the index of a twisted curve
at a node is invariant under specialization, it is enough to consider the case
where $\cC$ has exactly one node. 

{\sc Case 1: a separating node.} Write $C = C_1 \cup C_2$, where $C_i$ are the
irreducible components, and let $p_i\in C_i$ be the points over the node.
Consider the 
restriction of $P \to \cC$ to the general locus $C\gen$. Since 
every tame abelian \'etale cover of $C_i -p_i$ extends to an \'etale cover of
$C_i$, we have that $ P \to C$ is unramified, therefore $\cC = C$.

{\sc Case 2: a nonseparating node.}  This case is especially simple in
case $m$ is prime: we have $H^1 (C, \ZZ/m\ZZ)=(\ZZ/m\ZZ)^{2g-1}$, thus
$C$ admits no  connected $G$-cover. 
Since  $\cC$ does have a connected $G$-cover, it is not isomorphic to
$C$, and therefore the node must be twisted. Since the order of the
stabilizer at the 
node divides the exponent $m$ of $G$, the assumption that $m$ is prime
implies that this order is
precisely $m$.

In general,
consider the Leray spectral sequence 
of \'etale cohomology groups for 
$\pi: \cC \to C$: 

$$ H^i(C, \bR^j\pi_*( \ZZ/m\ZZ)) \Longrightarrow H^{i+j} (\cC, \ZZ/m\ZZ)$$

Note that $\pi_* (\ZZ/m\ZZ) = \ZZ/m\ZZ$, therefore we have an exact
sequence

$$ 0\lrar H^1(C, \ZZ/m\ZZ) \lrar  H^1(\cC, \ZZ/m\ZZ) \lrar H^0(C,
\bR^1\pi_* (\ZZ/m\ZZ)) $$

The existence of the twisted $G$-cover $P \to \cC$ shows that 
$G=(\ZZ/m\ZZ)^{2g} \subset H^1(\cC, \ZZ/m\ZZ)$, therefore the order of $
H^1(\cC, \ZZ/m\ZZ)$ is at least $m^{2g}$.

We have  $H^1(C, \ZZ/m\ZZ)=(\ZZ/m\ZZ)^{2g-1}$. By Proposition
\ref{Prop:hdi-to-moduli} in the appendix,
$\bR^1\pi_*(\ZZ/m\ZZ)$ is a 
sheaf concentrated at the node, whose stalk is $H^1(\Gamma, \ZZ/m\ZZ)$, where
$\Gamma$ is the stabilizer of a geometric point of $\cC$ over the node. Note
that $\Gamma$ is a cyclic group of order $m'$ dividing the exponent of
$G$, namely 
$m$. So the order of $H^1(\Gamma, \ZZ/m\ZZ)$ is precisely $m'$. The
exact sequence above becomes 
$$ 0\lrar (\ZZ/m\ZZ)^{2g-1} \lrar H^1(\cC, \ZZ/m\ZZ) \lrar \ZZ/m'\ZZ,
$$ therefore the order of $H^1(\cC, \ZZ/m\ZZ)$ is at most  $m^{2g-1}m'$. 
Combining the two inequalities,
we have that the order $m'$ of $\Gamma$ is precisely $m$. \qed


\subsection{The local system}
\begin{proposition} \label{Prop:local-system} 
 Given a pre-level-$m$ curve  $h:\cC \to S$, the sheaf 
$\bR^1h_*\ZZ/m\ZZ$ is a local system.
\end{proposition}

\begin{lemma}\label{Lem:tei-exists}
 A pre-level-$m$ curve over a strictly henselian ring
 has a twisted Teichm\"uller $\level$-structure.
\end{lemma}

\proof Let $\cC \to S$ be a
pre-level-$m$ curve over the spectrum of a strictly henselian ring, with closed
fiber $\cC_0$.  By the Proper-Base-Change Theorem for tame stacks
(\ref{Th:proper-base-change} in the appendix) we have $H^1(\cC, G) =
H^1(\cC_0, G)$, implying that any 
twisted $G$-cover on $\cC_0$ extends to a  twisted $G$-cover, inducing a 
 twisted Teichm\"uller structure on $\cC$. Therefore we may assume that $S$ 
is the spectrum of an algebraically closed field.

We  note that a pre-level-$m$ curve $\cC$ over an algebraically
closed 
field $k$ is determined up to isomorphisms by the underlying
curve. This follows since given a 
chart $(U, \Gamma)$ for $\cC$ at a nonseparating node, the strict
henselization of  $U$ is 
the 
strict henselization of $\Spec k[\xi, \eta]/(\xi\eta)$ and the action  of
some generator of $\Gamma$ is given by $(\xi, \eta) \mapsto
(\zeta_m\xi, \zeta_m^{-1} \eta)$.

Given a pre-level-$m$ curve $\cC_0$   over an algebraically
closed 
field, let $R$ be a discrete valuation ring with residue field $k$ and fraction
field $K$, and Let $\cC
\to \Spec R$ be a twisted curve deforming $\cC_0$, whose general fiber is
smooth. The geometric general fiber $\cC_{\bar\eta}$ admits a  Teichm\"uller
structure, therefore replacing $R$ by a finite extension we may assume
$\cC_\eta$ admits a  Teichm\"uller
structure. Since $\TCtei{G}$ is proper, replacing $R$ by a finite extension we
may assume $\cC_\eta$  has an extension $\cC'\to \Spec R$ with a Teichm\"uller
structure. Since the moduli spaces $C'$ and $C$  of $\cC'$ and $\cC$ are both
stable and have 
the same generic fiber, we have that $C'$ and $C$ are isomorphic. This implies
that the closed fiber of $\cC'$ is isomorphic to $\cC_0$ by the argument above.
This proves the Lemma. \qed

\marg{Little explanation on the morphism below is needed.}
\begin{lemma}\label{Lem:tei-isom}
Let $h:\cC\to S$ be a pre-level-$m$ curve with a twisted Teichm\"uller
$\level$-structure. The induced sheaf homomorphism $(\ZZ/m\ZZ)^{2g}_S \to
\bR^1h_*(\ZZ/m\ZZ)$ is an isomorphism.
\end{lemma}

\proof It is sufficient to check that this homomorphism induces isomorphisms on
stalks. By proper base change (Theorem \ref{Th:proper-base-change}) it
is enough 
to check this when 
$S$ is the spectrum of an algebraically closed field.

Let $\cC_0$ be such a pre-level-$m$ curve. Fix  a twisted Teichm\"uller
structure on $\cC_0$. It induces an injection $ (\ZZ/m\ZZ)^{2g} \hookrightarrow
H^1(\cC_0, \ZZ/m\ZZ)$.

Let $\cC' \to \Spec R$ be a deformation of $\cC_0$ with smooth generic
fiber, where $R$ a 
strictly 
henselian 
discrete valuation ring. Let $K$ be the fraction field and $K^*$ be its
separable closure. We denote by $i_0:\cC_0 \to \cC'$ the natural inclusion. We
have a cartesian diagram 
$$ \begin{array}{ccccc} 
 C^* & \lrar & C'_\eta & \lrar & \cC' \\
\dar &     &   \dar   &     & \dar  \\
\Spec K^*& \lrar & \Spec K & \lrar & \Spec R
\end{array} 
$$
Denote the composite morphism $\rho: C^* \to \cC'$. We have
$i_0^*(\rho_*\ZZ/m\ZZ) = \ZZ/m\ZZ_{\cC_0}$ (see \cite{SGA7}, Proposition
I.4.10).   
By proper base change (Theorem \ref{Th:proper-base-change}) we have
$$H^1 (\cC_0, 
\ZZ/m\ZZ) =H^1(\cC',  \rho_*\ZZ/m\ZZ), $$ and this injects in  the
group $H^1 (C^*, 
\ZZ/m\ZZ)$, which is isomorphic to  $ (\ZZ/m\ZZ)^{2g}$. This proves
the Lemma. \qed 

{\em Proof of Proposition \ref{Prop:local-system}.}  Let $\cC \to S$ be a
pre-level-$m$ curve. To show that   $\bR^1h_*\ZZ/m\ZZ$ is a local system we may
assume that $S$ is the spectrum of a strictly henselian ring. In this case the
statement follows immediately from Lemmas \ref{Lem:tei-exists} and
\ref{Lem:tei-isom}. \qed

The Proposition allows us to define a category  $\ocM_g^{(m)}$ of {\em twisted
curves with level $m$ structure} whose objects are 
pre-level-$m$ 
curves $h:\cC\to S$ along with  isomorphisms $(\ZZ/m\ZZ)^{2g}_S \to
R^1h_*\ZZ/m\ZZ$, and morphisms given by fibered squares as usual.

\begin{theorem}\label{Th:twisted-level-m} The category  $\ocM_g^{(m)}$
of twisted curves with level 
$m$ structure is an algebraic stack 
isomorphic to $\TCtei{\level}$.
\end{theorem}

\proof Given an object of $\TCtei{G}(S)$ over $\cC \to S$, we have by
Proposition 
\ref{Prop:pre-level} that $\cC \to S$ is a pre-level-$m$ curve. By Lemma
\ref{Lem:tei-isom} we have an isomorphism $(\ZZ/m\ZZ)^{2g}_S \to
R^1h_*\ZZ/m\ZZ$, giving an object of  $\ocM_g^{(m)}(S)$. 

In the other direction, let $\cC \to S$ be a pre-level-$m$ curve, along with an
isomorphism $(\ZZ/m\ZZ)^{2g}_S \to R^1h_*\ZZ/m\ZZ$.  There is an \'etale
surjective map $S' \to S$ such that this isomorphism comes from a group
homomorphism $(\ZZ/m\ZZ)^{2g} \to H^1(\cC',\ZZ/m\ZZ)$, where $\cC' = \cC
\times_S S'$. This corresponds to a principal $G$-bundle $P \to \cC'$.
Denote $S'' = S'\times_S S'$ and $\cC'' \to S''$ the pullback. The two
pullbacks $P_i\to \cC''$ of $P$ become isomorphic on some \'etale cover $T''
\to S''$. An isomorphism over $T''$ descends to an isomorphism of $P_1$ with
$P_2$ {\em in the catergory of twisted Teichm\"uller structures $\TCtei{G}$},
giving descent data for $P \to \cC'$ to an object of  $\TCtei{G}(S)$.
 \marg{Need to do arrows.}
\qed

\section{Automorphisms of twisted $G$-covers} 

We start this section with a  concrete description of the group of
automorphisms of a twisted 
curve $\cC$ acting trivially on the coarse curve $C$. We then turn to
automorphisms of $G$-covers, and show that, in case $G$ surjects to
$\level$, $m\geq 3$, every $G$-automorphism of a $G$-cover acts trivially on the
coarse curve $C$. We give some structure results on this automorphism group in
case $G$ is a characteristic quotient, and construct some fine
moduli spaces of twisted Teichm\"uller $G$-structure and of twisted
$G$-covers. 

\subsection{Automorphisms of twisted curves}
Let $\cC$ be a twisted curve over an algebraically closed field, and let $C$ be
its moduli space. For each node $x\in C\sing$ denote by $\Gamma_x$ the
stabilizer of a geometric point of $\cC$ over $x$, which is a cyclic group.

\begin{proposition}\label{Prop:aut-twisted-curve}
Denote by $\Aut_C(\cC)$ the automorphism group of $\cC$ over $C$ in
the category of 
twisted curves.
There is an isomorphism $$\Aut_C(\cC)\simeq \prod_{x\in C\sing}\Gamma_x.$$
\end{proposition}

In other words, every node contributes exactly $\Gamma_x$ to this
automorphism group (and the markings do not).

We use the following lemma, whose statement and proof are identical to
the classical case:
\begin{lemma}\label{Lem:aut-conn-etale-cover}
Let $P\to \cC$ be a connected \'etale $G$-cover of a twisted
curve. Then $\Aut_\cC P = G$ 
and $\Aut_\cC^G P = Z(G)$. 
\end{lemma}
{\em Proof of the Lemma.}\marg{Need to clarify $G$ vs. $G^{opp}$}
Clearly $\Aut_\cC P \supset G$. An 
automorphism $\phi\in \Aut_\cC P$ (which we write acting on the right)
pulls back to an automorphism of the
trivial cover $P\times_\cC P  \to P$ commuting with the Galois group
of the base change $P \to \cC$, which contains $G$ (which, to keep
things compatible, acts on the left). Restricting to a geometric
point and identifying the fiber with $G$, we get an element of the
permutation group of the set $G$ commuting with the action of $G$ on
the left, which therefore is an element of $G$ acting on the right,
giving the first claim.

If, moreover, $\phi$ commutes with the action of $G$ on the right,
then it is an element of the center $Z(G)$, giving the second claim.
\qed 

{\em Proof of the Proposition.} We can view $\Aut_C\cC$ as the group
of global sections of the  
\'etale sheaf of relative automorphisms - see Lemma \ref{Lem:aut-sheaf} below. 
Over $C\gen$ the map $\cC \to C$ is 
an isomorphism, and the sheaf is trivial. Therefore the automorphism
group is a product of 
contributions form small \'etale neighborhoods of twisted nodes and markings.
Focusing on one of these twisted points $x$, we may replace $\cC$ by an
affine twisted 
curve having 
the same local picture at a twisted geometric point, that is, $\cC =
[U/\Gamma_x]$ where $U$ is either $\Spec k[z]$ or $\Spec k[z,w]/(zw)$,
with $\Gamma_x = \bmu_r$ acting as described in
\ref{Sec:at-marking} or  \ref{Sec:at-node}, respectively.

Consider the exact sequence
\begin{equation}\label{exact-sequence} 1 \lrar \Aut_{[U/\Gamma_x]}U
\lrar \Aut_C (U \to 
[U/\Gamma_x]) \lrar 
\Aut_C\cC.
\end{equation}

\begin{claim} We have $\Aut_{[U/\Gamma_x]}U = \Gamma_x.$ 
\end{claim}

{\em Proof of claim.} 
This follows from  Lemma \ref{Lem:aut-conn-etale-cover},
 since $U$ is a connected \'etale $\Gamma_x$-cover of
${[U/\Gamma_x]}$.
\qed

\begin{claim} The canonical inclusion $\Aut_C (U \to [U/\Gamma_x]) \subset
\Aut_C U $ is an isomorphism.
\end{claim}

{\em Proof of claim.} 
In case $U = \Spec k[z]$, we have $\Aut_C (U \to [U/\Gamma_x]) \subset
\Aut_C U = \Gamma_x$. The exaxt sequence (\ref{exact-sequence})
implies that equality holds.  

In case $U=\Spec k[z,w]/(zw)$, we have $\Aut_C (U \to [U/\Gamma_x]) \subset
\Aut_C U = \bmu_r^2$, where the action of $(\zeta_1,\zeta_2)\in
\bmu_r^2$ is via $$ (z,w)\ \  \mapsto\ \  (\ \zeta_1\ z\ ,\ \zeta_2\ w\ ).$$
This action clearly commutes with the action of $\Gamma_x$, which
means that $(\zeta_1,\zeta_2)$ acts on $U \to [U/\Gamma_x]$. \qed

\begin{claim} The morphism on the right in the sequence
(\ref{exact-sequence})  is   
surjective.  
\end{claim}

{\em Proof of Proposition assuming the claim.} 
\begin{itemize}
\item In case $U = \Spec k[z]$, 
$\Aut_C\cC = \Gamma_x/\Gamma_x$ is trivial.    

\item In case $U=\Spec k[z,w]/(zw)$, we have $\Aut_C\cC = \bmu_r^2/\Gamma_x
\simeq \Gamma_x$, as required.\qed
\end{itemize}


{\em Proof of the claim.}  Let $\phi\in \Aut_C[U/\Gamma_x]$. This
comes from a functor 
$[U/\Gamma_x] \to 
[U/\Gamma_x]$ preserving $U/\Gamma_x$. 

The canonical morphism $U \to [U/\Gamma_x]$ corresponds to the diagram 
$$ \begin{array}{ccc} 
U \times \Gamma_x & \lrar & U \\
\dar &&\\
U \end{array} $$
where the vertical map is the projection on the first factor and the horizontal
map is the action of $\Gamma_x$ on $U$. The automorphism $\phi_U$ gives another
principal $\Gamma_x$ bundle 
$$ \begin{array}{ccc} 
P & \lrar & U \\
\dar &&\\
U \end{array} $$

Since $U$ has trivial tame fundamental group,  we may choose a section $U \to
P$, and composing with the horizontal map 
$P \to U$ 
we  obtain an automorphism of $U$ over $U/\Gamma_x$. This gives a
lifting of $\phi$ to $\Aut_CU$, which by the previous claim is the
same as $\Aut_C(U \to \cC)$. \qed






%

In the proof we used the following lemma:

\begin{lemma} \label{Lem:aut-sheaf}
Let $\cX$ be a separated Deligne--Mumford stack  over a scheme
$S$. Suppose that there is an open and scheme-theoretically dense
substack of $\cX$ which is an algebraic space. Then the functor that sends
each \'etale morphism of finite type $U \to S$ to the group of
automorphisms $\Aut_U(\cX_U)$ is a sheaf on the small \'etale
site of $S$.
\end{lemma}

The group $\Aut_U(\cX_U)$ is the group of base-preserving equivalences of
categories of $\cX_U$ with itself, modulo isomorphism. Recall that the
groupoid of base-preserving equivalences of categories of $\cX_U$ with
itself is in fact equivalent to a group (\cite{stablemaps}, Lemma
4.2.3). In other words, no such equivalence has nontrivial automorphisms.

\proof First of all, let us check that if $\{U_i \to U\}$ is an \'etale
cover, $F, G \colon \cX_U \to \cX_U$ are base-preserving equivalences,
such that their pullbacks $F_{U_i}, G_{U_i} \colon \cX_{U_i} \to
\cX_{U_i}$ are all isomorphic, then $F$ and $G$ are isomorphic.

This follows immediately from the fact that the isomorphisms of pullbacks
of $F$ and $G$ to \'etale morphisms into $U$ form a sheaf in the small
\'etale topology of $U$. In concrete terms, given isomorphisms $\phi_i \colon
F_{U_i} \simeq G_{U_i}$, the pullbacks of $\phi_i$ and $\phi_j$ to
$F_{U_{ij}} \simeq G_{U_{ij}}$ must coincide, because of the unicity of
isomorphisms, therefore $\{\phi_i\}$ satisfies the cocycle condition,
and, by the stack axioms,
the $\phi_i$ descend to an isomorphism of $F$ with $G$.

Now assume that you are given a collection $F_i \colon \cX_{U_i} \to
\cX_{U_i}$ of base preserving equivalences, such that $(F_i)_{U_{ij}}$ and
$(G_i)_{U_{ij}}$ are isomorphic. Let $\phi_{ij} \colon (F_j)_{U_{ij}}
\simeq (F_i)_{U_{ij}}$ be an isomorphism; the unicity of isomorphisms
insures that the cocycle condition
$$
\phi_{ij} \phi_{jk} = \phi_{ik} \colon (F_i)_{U_{ij}} \simeq
(F_i)_{U_{ij}}
$$
is satisfied.

Let $T \to U$ be a morphism, and set $T_i = T \times_U U_i$ and $T_{ij} =
T \times_U U_{ij}$. Suppose that $\xi$ is an object of $\cX(T)$. Since we
have $F_i(\xi_{T_i})_{T_{ij}} = F_i(\xi_{T_{ij}})$, the isomorphism
$\phi_{ij}(\xi_{T_{ij}})
\colon F_j(\xi_{T_{ij}}) \simeq F_i(\xi_{T_{ij}})$ yield isomorphisms
$\psi_{ij}\colon F_j(\xi_{T_j})_{T_{ij}} \simeq F_i(\xi_i)_{T_{ij}}$. The
cocycle
condition on the $\phi_{ij}$ says that these isomorphisms $\psi_{ij}$
give descent data; therefore we obtain an object $F(\xi)$ of $\cX(T)$,
together with isomorphisms $F(\xi)_{U_i} \simeq F_i(\xi_{U_i})$.

It is a simple matter to check that if $f \colon\xi \to \eta$ is an arrow
in $\cX(U)$, then the restrictions $f_i \colon F_i(\xi_{U_i}) \to
F_i(\eta_{U_i})$ glue together to yield an arrow $F(f) \colon F(\xi) \to
F(\eta)$; therefore we obtaine a functor $F \colon \cX_U \to \cX_U$, which
lifts to the $F_i$, as desired.
 \qed

\subsubsection{Automorphisms and
deformations}\label{Sec:automorphisms-and-deformations}   In case $\cC$ is
balanced and unmarked, one can read off the automorphism group from
the deformation space. Again, the problem is local, so we may focus on
the case where $\cC$ has a unique twisted node. We note
that the deformation space $\Def_\cC$ of $\cC$ surjects to the
deformation space 
$\Def_C$ of $C$. Also, since there is a unique isomorphism class of
twisted curves  with given twisting having coarse moduli space $C$, we
have that $\Def_C = \Def_\cC / \Aut_C\cC$. Moreover, since the generic
curve is smooth and since there are no markings, it is untwisted. This
means that the action of $\Aut_C\cC$ on $\Def_\cC$ is effective. 

From the analysis of  the sequence (\ref{deformation-sequence}) in section
\ref{Sec:deformation-theory} we see that we only need compare the
deformation spaces of the nodes on $\cC$ and $C$. Clearly the map from
the deformation space $\Spec k[[t]]$ of $U=\{zw=0\}$ (having versal
family $zw=t$)  to the
deformation space $\Spec k[[s]]$ of $U/\mu_r = \{xy=0\}$ (having versal
family $xy=s$) is given by $s = t^r$, and its Galois group is
$\bmu_r = \Gamma_x$. 

\subsection{Serre's lemma and existence of tautological families} In the rest 
of this section we
 assume $n=0$, that is, the curves have no markings.  We fix $G$, a
 finite group, and write $\bbS = 
\Spec \ZZ[1/\# G]$.

The following is an interpretation of a well known lemma of Serre:

\begin{lemma} Assume $G$ admits a surjection onto $(\ZZ/m\ZZ)^{2g}$, for some
$m\geq 3$. Then every $G$-automorphism of a balanced, connected
twisted $G$-cover $P 
\to \cC$  acts trivially on the coarse curve $C$. In other words,
$$\Aut^G(P \to \cC) = \Aut^G_C(P \to \cC).$$
\end{lemma}
\proof Let $Q \to \cC'$ be the twisted $\level$-cover associated to
$P\to \cC$ obtained using the functoriality result
\cite{stablemaps}, Corollary 9.1.2, applied with $\cM = \cB G$ and
$\cM'=\cB \level$.  An automorphism of the
twisted $G$-cover $P \to \cC$ gives an 
automorphism  of $Q \to \cC'$ over an automorphism $\phi:\cC' \to \cC'$. Recall
that $Q \to \cC'$ gives a basis for $H^1(\cC', \ZZ/m\ZZ)$ (see Lemma
\ref{Lem:tei-isom}). Since $\phi^*Q \simeq Q$, we have that $\phi^*: H^1(\cC',
\ZZ/m\ZZ)\to H^1(\cC', \ZZ/m\ZZ)$ is the identity. But $H^1(C,
\ZZ/m\ZZ)\subset H^1(\cC', \ZZ/m\ZZ)$, therefore $\phi$ induces the identity on
$H^1(C,\ZZ/m\ZZ)$. By Serre's Lemma for stable curves (see, e.g.,
\cite{A-O:alterations}, Lemma 3.5.) it follows that $\phi$ induces the identity
on $C$. \qed

We deduce the following well known corollary:
  
\begin{corollary}  Assume $G$ admits a surjection onto
$(\ZZ/m\ZZ)^{2g}$, for some $m\geq 3$. Then the
morphism of coarse moduli spaces $\bB^{\tei}_g(G)\to \obM_g$ admits a
lifting $\bB^{\tei}_g(G)\to \ocM_g$. In other words, $\bB^{\tei}_g(G)$
carries a tautological family of stable curves. 
\end{corollary}
\proof Let $\cC_0$ be a twisted curve over an algebraically closed
field admitting a twisted 
Teichm\"uller $G$-structure corresponding to a connected balanced
$G$-cover $P \to 
\cC_0$. Denote $A^G = \Aut^G(P \to \cC_0)$. The local picture of
$\cB^{\tei}_g(G)$ at   
the point corresponding to $P \to \cC_0$ is $[\Def_{P\to \cC_0} /
A^G]$, and the local picture of the universal curve is $[\cC/ A^G]$,
where $\cC \to \Def_{P\to \cC_0}$ is the twisted curve 
underlying the universal
deformation. Consequently, on the level of coarse moduli spaces, we 
have that the local picture of $\bB^{\tei}_g(G)$ at  
the point corresponding to $P \to \cC$ is the scheme $\Def_{P\to \cC} /
A^G$, and the local picture of the coarse moduli space of the
universal curve is 
$C/ A^G$, where $C$ is the
{\em coarse} curve underlying the universal deformation. Since the
action of $A^G$ on 
$C_0$ is trivial, we have that $C/ A^G \to  \Def_{P\to \cC} /
A^G$ is a stable curve, as required. \qed

\subsection{Structure of automorphisms of connected twisted $G$
covers} 
Let $P \to \cC$ be a  connected twisted $G$ cover over
an algebraically closed 
field. Given a node $x$ of $C$, we denote by $r_x$ the index of $\cC$ at
$x$. 

We wish to have some understanding the $G$-automorphism group
$\Aut_C^G(P \to \cC)$ of the twisted 
$G$-cover $P \to \cC$ fixing $C$. It is easy to see that $\Aut_C^G(P
\to \cC) = \Aut_C^G P$, since $\cC$ can be recovered as $[P/G]$. One
may try to study it via its natural 
embedding  as the centralizer of $G$ in the group $\Aut_C P$, but the
latter group is in general too big - the action of an element of
$\Aut_C P$ is in general not compatible with local charts for the
twisted cover $P \to \cC$. 

We denote $A = \Aut_C (P\to \cC)$, the automorphism group of the morphism
$P\to \cC$,  
fixing $C$ (but not necessarily commuting with $G$). This is precisely
the group of automorphisms of $P$ over 
$C$ preserving the charts of $P \to \cC$ as a twisted cover.
  Then the $G$-equivariant automorphisms are 
$A^G$. These in turn can be thought of as  the
$G$-automorphisms of $P$ over $C$.

Set 
$M 
\stackrel{\operatorname{def}}{=} 
\displaystyle \prod_{x \in C\sing} \Gamma_x $. 
 
\begin{lemma} We have an exact sequence 
$$ 1 \lrar G \lrar A \lrar M $$
which, when taking $G$-invariants, gives an exact sequence
$$ 1 \lrar Z(G) \lrar A^G \lrar M $$
\end{lemma}

\begin{proof}
Consider the natural sequence
$$ 1 \lrar \Aut_\cC P \lrar \Aut_C (P\to \cC) \lrar \Aut_C \cC.$$
The group 
$\Aut_\cC P$ is naturally isomorphic to $G$ since $P$ is connected and $P\to
\cC$ is a principal bundle (Lemma \ref{Lem:aut-conn-etale-cover}). The term
on the right is $\prod_{x \in C\sing} 
\Gamma_x $ by Proposition \ref{Prop:aut-twisted-curve}.
\end{proof}

Recall that the stack of balanced twisted $G$-covers is nonsingular and flat
over $\ocM_g$. This implies that any object  $P \to \cC$  can be
deformed to a 
smooth object in characteristic 0. It also implies that $P \to \cC$
can be deformed to an object in characteristic 0 preserving the topological
type of $\cC$.

\begin{definition} We say that a connected balanced twisted $G$-cover
$P \to \cC$ is {\em 
characteristic}  
if for some deformation to a smooth $G$ cover  $P' \to C'$ in characteristic 0,
and some choice of base point $s$ in $C'$, the kernel of the
corresponding epimorphism $\pi_1^{geom}(C',s) \to G$ is 
a characteristic subgroup of $\pi_1^{geom}(C',s)$.
\end{definition}

 Since the stack
$\TC{g}{\bal}{G}$ is smooth over $\bbS$, and since the property of
$\pi_1^{geom}(C',s) \to G$ being characteristic is Galois invariant, this 
property is independent of the choice of  deformation, and thus it is an
invariant of the connected component of the stack.

\begin{lemma} \label{Lem:char-dehn}
Assume that $P \to \cC$ is characteristic. Then the homomorphism $ A \lrar M $
is surjective, giving an exact sequence 
$$ 1 \lrar G \lrar A \lrar M \lrar 1. $$ 
Furtheremore, when the base field is $\CC$, then for any deformation 
 $\cC_\Delta\to \Delta$ of $\cC$ with smooth generic fiber, and each
node $x\in C$, there 
is a  generator $\sigma_x\in \Gamma_x\subset M$ with a lifting
$\delta_x\in A$, whose action on 
$G$ via conjugation is obtained by the action of a Dehn twist $D_x$ along the
vanishing cycle of the node, on the fundamental group of a nearby
smooth curve in the deformation $\cC_\Delta\to \Delta$.
\end{lemma}
We remark that this exact sequence can be shown to  split. 
Our proof relies on topological considerations over $\CC$. 
\subsubsection{Reduction to $\CC$} We claim that it suffices to show
the Lemma when $k = 
\CC$. First consider the case of characteristic 0. Note that $A$,
being the group of points of a finite group scheme, is invariant under
extensions of algebraically closed field. It follows that if the
statement holds in over $\CC$ 
then it holds in characteristic 0. 

If $P_0 \to \cC_0$ is in positive
characteristic, we choose a deformation with constant topological type $P \to
\cC$ on a discrete valuation ring $R$ of mixed characteristic, with fraction
field  $K$. We may assume $R$ contains all the roots of unity of order
dividing the order of $G$. The sequence above comes from a sequence of
group schemes $ 1 \lrar 
\tilde G \lrar \tilde A \lrar \tilde M $, Here $\tilde G$ is assumed 
constant, and $\tilde M = \prod \bmu_{r_x}$  is constant since the
deformation has constant 
topological type and $R$ contains the roots of unity of order
$r_x$. An element 
$m\in M$ is therefore also an element of $\tilde M(K)$. This element
lifts to an 
element of $\tilde A(\bar K)$, since $K$ has characteristic 0, which
is the case
we settled above. This     
implies that after replacing $R$ by a finite extension, the element $m$ lifts
to $\tilde A(K)$. Since $\tilde A$ is a finite group scheme, the
element in $\tilde A(K)$ specializes
to an element $a$ of $A$.

\subsubsection{Topologically trivial parametrization}
Now consider the case   $k = \CC$. Let $P_\Delta \to \cC_\Delta \to \Delta$ be
a small 
analytic deformation of $P\to \cC$ with smooth generic fiber, where
$\Delta$ is the disc around the 
origin in $\CC$. We assume that the only singular fiber lies over the
origin in $\Delta$.  Let $\cC_{1}$ be the fiber over a general point $t_1\in
\Delta^*$. 
Consider the parametrized line segment $\beta:[0,1] \to \Delta$ given
by $\beta(t) = t t_1$,
connecting $0$ with 
$t_1$. 

\begin{claim}
There is a family of continuous maps $\psi_t:\cC_1  \to \cC_{\beta(t)}$, such
that $\psi_1$ is the identity, $\psi_t$ is a homeomorphism whenever $t\neq 0$,
and $\pi\circ\psi_0:\cC_1 \to C_0 = C$ is the contraction of the vanishing
cycles.
\end{claim}
{\em Proof of claim.} First recall a typical construction of a
continuous map  $\phi_t:\cC_1  \to C_{\beta(t)}$ to the {\em coarse}
curve, using polar coordinates. Write $C_\Delta = C' \cup \bigcup V_i$, where
$C'$ is an open set which 
is  topologically
trivial over $\Delta$, and $V_i$ are small neighborhoods of the nodes of
$C_0$. We may assume that we have an analytic isomorphism $V_i \simeq
\left\{(x_i,y_i, u) : |x_i|, |y_i| \leq 1, u\in \Delta, x_iy_i
= h_i(u) 
\right\}$, for some analytic 
function $h_i$ on $\Delta$. We note that the existence of $\cC_\Delta$
implies that $h_i = g_i^{r_i}$.

We focus on one of the open sets $V_i$ and drop the subscripts $i$ for
simplicity of notation.

For $t \in (0,1]$ write $h(\beta(t)) = \tau(t) e^{2\pi i \alpha(t)}$,
where $\tau(t)\in \RR_{\geq 0}$ and $\alpha(t) \in \RR/\ZZ$. Since
$\beta$ is 
linear, this extends continuously to $t = 0$ as well. 

We use coordinates $(\xi,\eta,\theta)$ on $\RR_{\geq 0}^2\times
\RR/\ZZ$. Consider $\tilde V\subset \RR_{\geq 0}^2\times \RR/\ZZ$, the
inverse image  
of $V$ by the polar coordinates map  
$$\begin{array}{ccl} x & = & \xi\ e^{\,2\pi i\,\theta}\\
                     y & = & \eta\ e^{\, 2\pi i\, (\alpha(t)\ -\ \theta)}.
\end{array}
$$
The inverse image of the point $x=y=0$ is a circle representing the
vanishing cycle. 

The fibration $\tilde V \to [0,1]$ is  topologically
trivial: all the fibers are homeomorphic to a cylinder.

Gluing each $\tilde V_i$ into $\beta^*(C_\Delta)$ instead of $V_i$, we
get a topologically 
trivial fibration $C^{polar} \to [0,1]$, and choosing a trivialization
we have a continuous map 
$\cC_1\times [0,1] \to C_\Delta$.

We now lift this map to the {\em stack} $\cC_\Delta$. Locally near a
node, $\cC_\Delta$ is 
given by $U: \{zw = g(t)\}$, where $x = z^r, y=w^r$ and $h(t) = g(t)^r$. Over
$\tilde V$ we have a corresponding {\em \'etale}  covering given by
$(\zeta,\omega,\gamma)$, where
$$\begin{array}{ccl} \xi & = & \zeta^r\\
                     \eta & = & \omega^r\\
                     \theta & = & r\cdot \gamma
\end{array} 
$$
which maps down to $U$ via
$$\begin{array}{ccl} z & = & \zeta\ e^{\,2\pi i\,\gamma}\\
                     w & = & \omega\ e^{\, 2\pi i\, (\alpha(t)/r\ -\ \gamma)}.
\end{array}
$$
and this map is clearly $\bmu_r$-equivariant. This gives the desired
lifting $\psi_t:\cC_1 \to \cC_\Delta$. \qed

It should be remarked that, using polar coordinates on both $\Delta$
and $C_\Delta$, one can define a topologically locally trivial
fibration $C^{polar}_{\Delta} \to \Delta^{polar}$. This is a morphism
of real-analytic manifolds with corners, which is an instance of the
logarithmic space associated to the natural logarithmically smooth
structure on 
$C_\Delta \to \Delta$, see \cite{Kato-Nakayama}\marg{Precise citation needed}.

We continue the proof of  the lemma.
The pullback by $\psi_t$ of $P_t\to \cC_t$ gives a topological principal bundle
over $\cC_1 \times [0,1]$, which therefore must be constant. This means that
$\psi_0^* P_0$ is isomorphic as a topological $G$-bundle to $P_1$, and
therefore is characteristic.   

\subsubsection{Making space around the vanishing cycle} Since Dehn
twists are not analytic in nature, and since we want to keep some
analytic properties, we make some space around every vanishing cycle
in $\cC_1$ by inserting a cylinder, inside which all the non-analytic
activities will occur.

Fix a node $x\in C\sing$.  A small neighborhood $W_x$ of the vanishing cycle
$(\pi\circ\psi_0)^{-1}x$ is homeomorphic to an open cylinder
$S^1\times (-\epsilon, \epsilon)$. We 
replace this by 
$$W_x^{\cyl} = S^1 \times \bigl(\  (-\epsilon,0] \cup [0,1] \cup
[1, 1+\epsilon) \ \bigr) = S^1\times (-\epsilon, 1+ \epsilon).$$
 Doing this at each node, we obtain a topological
surface 
$C^{\cyl}$ with a continuous map $\eta:C^{\cyl}\to \cC_1$ shrinking
each 
$S^1 \times 
[0,1]$ to the cycle $S^1\times \{0\}$. Denote $\psi_0\circ \eta = \psi^{\cyl}:
C^{\cyl}\to \cC_0$ and let $P^{\cyl}=(\psi^{\cyl})^* P_0$. This is a
characteristic topological cover, in the sense that the corresponding subgroup
$\pi_1(P^{\cyl}) \subset \pi_1(C^{\cyl})$ is characteristic. 

\subsubsection{The Dehn twist} Fix one node $x \in C$.
Define the following homeomorphism $D_x:C^{\cyl}\to C^{\cyl}$. On the
complement 
of  $S^1 \times [0,1] \subset W_x^{\cyl}$ it is defined to be the
identity. On  
$S^1 \times [0,1]$ it is defined by $D_x(z,t) = (e^{\,2\pi i\, t} z, t)$. The
induced action of $D_x$ on $\pi_1 (C^{\cyl}) \simeq \pi_1 (\cC_1)$ is precisely
the Dehn twist associated to the vanishing cycle over $x$.

Consider the pullback $D_x^*P^{\cyl}$. Since $P^{\cyl}$ is a characteristic
topological cover, this pullback is isomorphic to $P^{\cyl}$ as a topological
cover, so $D_x$ lifts to a homeomorphism  $\delta^{\cyl}_x:P^{\cyl}\to
P^{\cyl}$.

\subsubsection{The action on $\cC$}
We claim that $\delta^{\cyl}_x$ descends to an automorphism $\delta_x:P\to P$
whose image in $M$ is a generator of $\Gamma_x$.

First we claim  that there is a commutative diagram of topological stacks
$$\begin{array}{lcl}
C^{\cyl} & \stackrel{\psi^{\cyl}}{\lrar} & \cC \\
 D_x \dar&                            & \dar \sigma_x \\
C^{\cyl} & \stackrel{\psi^{\cyl}}{\lrar} & \cC
\end{array} $$ 
for a suitable  automorphism $\sigma_x$ corresponding to a
generator of $\Gamma_x\subset \Aut_C\cC$. To see this, consider a chart $(U,
\Gamma_x)$ for $\cC$ at $x$ with a lifting $W_x^{\cyl} \to [U/\Gamma]$ of
$\psi^{\cyl}:C^{\cyl} \to \cC$, and let 
$$\begin{array}{ccc}
\widetilde{W^{\cyl}_x} & \lrar & U \\
\dar & &\dar\\
W^{\cyl}_x & \lrar & [U/\Gamma_x]
\end{array}$$
be the cartesian diagram. Then $D_x$ lifts to an automorphism of
the cylinder $\widetilde{W^{\cyl}_x}$ which is trivial on $S^1 \times
(-\epsilon, 
0]$, acts on $S^1 \times [0,
1]$ by $(w, t) \mapsto (e^{\,2\pi i\, t/r_x} w, t)$, and rotates $S^1
\times [1, 1+\epsilon]$ by $e^{2\pi i / r_x}$ where $r_x$ 
is the order of $\Gamma_x$.

It follows that $\psi^{\cyl}\circ D_x = \sigma_x 
\circ \psi^{\cyl}$, where 
$\sigma_x\in \Aut_C\cC$ acts as the identity on the branch of $U$
under   $S^1 \times 
(-\epsilon, 
0]$  and as $e^{2\pi i / r_x}$ on the branch of $U$ under $S^1 \times [1,
1+\epsilon)$. 

\subsubsection{The action on $P$}
Now it is clear that $\delta^{\cyl}_x$ descends to a homeomorphism
$\delta_x:P\to P$. This homeomorphism is analytic outside of the nodes, and by
continuity it is analytic everywhere. Now we have a commutative diagram
$$\begin{array}{ccc}
P & \stackrel{\delta_x}{\lrar} & P \\
\dar & & \dar \\
\cC& \stackrel{\sigma_x}{\lrar} & \cC
\end{array}$$
which is what was required. \qed

\subsubsection{An algebraic proof}

Here is a sketch of an algebraic proof proposed by Johan de Jong. It is presented
geometrically, but one can easily make it algebraic, for instance
using the methods of \cite{Pikaart-DeJong}. 

Let $P_0\to \cC_0$ be a connected balanced stable twisted $G$-cover over $\CC$,
and let $C_0$ be the coarse stable curve. Let $\oC \to\oDelta$ be the
universal deformation of $C_0$, and let $P \to \cC \to \Delta$ be the
universal deformation of $P_0\to \cC_0$. We denote by $C$ the coarse moduli
space of $\cC$.

The deformation $C \to \Delta$ of $C_0$ induces a morphism $\Delta\to
\oDelta$. By \ref{Sec:automorphisms-and-deformations}  we have that
$\Delta\to \oDelta$ is Galois with 
Galois group $M=\Aut_{C_0}\cC_0$, where $M$ acts on $\cC \to \Delta$
by the universal property of the deformation.

Since $C = \oC\times_\oDelta \Delta$, the action of $M$ on $\Delta$
lifts as a product action to $C$. This is nothing but the action
induced by the action of 
$M$ on $\cC$.

Fix $\sigma\in M$ and let $t$ be a geometric point of $\Delta$ in the
smooth locus of $C \to \Delta$. The smooth curve $C_t$ carries two
covers $P_t \to C_t$ and $P_{\sigma (t)} = (\sigma^*P)_t \to C_t$. Since
the kernel of $\pi_1(P_t) \to \pi_1 (C_t)$ (with appropriate base
points) is a characteristic subgroup, these 
two covers are isomorphic, so there exists an isomorphism 
$$\delta_t:   P_t \to P_{\sigma (t)}$$
lying over 
$$\sigma:   C_t \to C_{\sigma (t)},$$ and, moreover, $\delta_t$ sends an
element of $G$ (acting on $P_t$) to an element of $G$ (acting on
$P_{\sigma (t)}$). 

Since the stack of stable curves of genus $g(P)$ is separated, this
reduces to an automorphism $\delta:P_0 \to P_0$ (lying over
$\sigma:C_0 \to C_0$, which is the identity). Moreover $\delta$ acts on
$G\subset \Aut_{C_0}P_0$, giving an automorphism of $[P_0/G]=\cC$
which clearly coincides with $\sigma$. Therefore the diagram 
$$\begin{array}{ccc} P_0 &\stackrel{\delta}{\lrar}  & P_0 \\
\dar && \dar \\
\cC_0 &\stackrel{\sigma}{\lrar}  & \cC_0\end{array}$$
is commutative, giving the surjectivity of $A \to M$.

  Denote by $\Delta^*, \oDelta^*$ the loci where the curves are
smooth, and by $P^* \to C^* \to \Delta$, respectively $\oC^* \to \oDelta$
the curve fibrations. On the level of fundamental groups, consider the
diagram with exact rows  and injective columns
$$\begin{array}{ccccccccc}
1& \to & \pi_1( P_t) & \to &\pi_1 (P^*) & \to & \pi_1 (\Delta^*) & \to&1\\
 &  & \dar &  &\dar & & \dar &  & \\  
1& \to & \pi_1( C_t) & \to &\pi_1 (\oC^*) & \to & \pi_1 (\oDelta^*) & \to &
1.\\ 
\end{array}$$
We have $M =  \pi_1 (\oDelta^*)\ / \ \pi_1 (\Delta^*) $, so an element
of $M$ lifts to an element $D\in \pi_1 (\oC^*)$, whose action on
$\pi_1(C_t)$ is the  Dehn twist (corresponding on the level of outer
automorphisms  to its image 
 in the 
monodromy group $\pi_1 (\oDelta^*)$). Since $\pi_1( P_t) \subset
\pi_1( 
C_t)$ is characteristic, we have $\pi_1( P_t) \subset \pi_1( \oC^*)$
normal, so $D$ acts on the quotient $G = \pi_1( C_t)\ / \ \pi_1( P_t)$
as a Dehn twist.

Incidentally, it can be shown that $\pi_1(P^*) \subset \pi_1(\oC^*)$
is normal with quotient group $A$.

\qed

Lemma \ref{Lem:char-dehn} implies that if $P\to \cC$ is
characteristic, there 
is a canonical 
induced homomorphism of groups from $M$ to the group 
$\Out (G)$ of outer automorphism of $G$. We can now summarize our
results as follows:

\begin{proposition} \label{Prp:group-is-trivial}
Assume $P\to \cC$ is a characteristic $G$-cover over an algebraically
closed field. Then
\begin{enumerate}
\item if the homomorphism $M \to \Out(G)$ is injective, then the subgroup
of the group of automorphisms of $P\to \cC$ as a twisted Teichm\"uller
structure, consisting of elements
acting trivially on $C$, is trivial.
\item If, moreover, $G$ admits a surjection to $(\ZZ/m\ZZ)^{2g}$ for
some $m\geq 3$, then the whole group
of automorphisms of $P\to \cC$ as a twisted Teichm\"uller structure
 is trivial.
\item If in addition the center $Z(G)$ is trivial then the group
of $G$-automorphisms $\Aut^G(P\to \cC)$ is trivial.
\end{enumerate} 
\end{proposition}

\proof The first statement follows from the exact sequence in Lemma
\ref{Lem:char-dehn} since the 
assumption implies that 
$A^G = Z(G)$. The second statement follows from Serre's Lemma. The
last statement follows from the same exact sequence since then $A^G$ is trivial.

\subsection{The groups of Looijenga and Pikaart--De Jong} 

We now describe two cases where the first two statements of the
proposition above hold. We note that condition (1) can be verified by
checking it for
twisted covers over $\CC$.  
We denote by  $\Pi_g$ be the fundamental group of a Riemann surface of genus
$g$.

\subsubsection{Looijenga's groups} Fix an integer $m\geq 3$. Choose a
Riemann surface $C$ of genus $g$, and let 
Let $C_{2} \to C$ be the maximal abelian \'etale cover of
exponent $2$. 
The Galois group $G_{2}$ of $C_{2}$ over $C$ is $H_1(C, 
\ZZ/2\ZZ)$, which is isomorphic to $(\ZZ/2\ZZ)^{2g}$. 
Let $C_{\prym{2m}{2}}\to C_{2}$ be the maximal abelian
\'etale cover of exponent $2m$.
The  curve $C_{\prym{2m}{2}}$ is a Galois cover of $C$. Denote its
Galois group by $G_{\prym{2m}{2}}$. This is clearly a characteristic
quotient of $\Pi_g$.

Let $C_0$ be a nodal curve and let $P_0 \to C_0$ be a connected
admissible $G_{\prym{2m}{2}}$-cover. We denote by $P_0 \to \cC_0$ the
corresponding twisted $G_{\prym{2m}{2}}$-cover. 
\begin{lemma}\begin{enumerate}
\item  If $x$ is a separating node of $\cC_0$ then
$\Gamma_x\simeq \bmu_{2m}$.
\item  If $x$ is a non-separating node of $\cC_0$ then
$\Gamma_x\simeq \bmu_{4m}$.\end{enumerate}
\end{lemma}
\proof
 There is an intermediate curve
$D_0$ which is a connected admissible $(\ZZ/2\ZZ)^{2g}$ cover of $C_0$.
According to \cite{Looijenga}, Proposition 2, every node of $D_0$ is a
non-separating node. If $x$ is a separating node of $\cC_0$, then
it follows from Proposition \ref{Prop:pre-level} that $D_0\to C_0$ is
 unramified at $x$, $P_0 \to D_0$ has ramification index $2m$,   and
 therefore $P_0\to C_0$ has ramification index $2m$. If $x$ is 
 non-separating, then it follows from the same proposition that
 $D_0\to C_0$ has  ramification index $2$ at $x$,  and $P_0\to D_0$
 has again index  $2m$, therefore $P_0\to C_0$ has index $4m$. The
 lemma follows. \qed 

Following Looijenga, let $E_0$ be the set of separating nodes of $C_0$
and $E_1$ the set of 
nonseparating nodes. Denote by $T= \prod_{x\in C\sing}\ZZ$ the group
of Dehn twists. According to Looijenga
\cite{Looijenga}, Proposition 3, the
kernel of the natural homomorphism $T \to \Out (G_{\prym{2m}{2}})$ is
precisely the  subgroup
$$T_0 = \prod_{x\in E_1} 2m\ZZ \times  \prod_{x\in E_2} 4m\ZZ.$$ It
follows from Lemma \ref{Lem:char-dehn} that we have an isomorphism 
$T/T_0 \simeq \prod_{x\in C\sing}\Gamma_x$. This implies that $M \to
\Out (G_{\prym{2m}{2}})$ is injective.

As an immediate outcome we have

\begin{theorem}
Suppose  $m\geq 3$. Then the moduli stack $\TCtei{G_{\prym{2m}{2}}}$
is a smooth 
projective scheme over $\ZZ[1/2m]$ admitting a finite flat morphism to
$\ocM_g$.  
\end{theorem}

\subsubsection{The groups of Pikaart-De Jong} 
 Following \cite{Pikaart-DeJong} we inductively define $\Pi_g^{(1)} =
\Pi_g$ 
and $\Pi_g^{(k+1)} = [\Pi_g^{(k)},\Pi_g]$, the group of $k$-th order
commutators. We denote by $\Pi_g^{(k),n} = \Pi_g^{(k)}\cdot \Pi_g^n$, where
$\Pi_g^n$ is the subgroup generated by $n$-th powers, and $G_g^{(k),n} =
\Pi_g/\Pi_g^{(k),n}$.  Since the exponent of $G_g^{(k),n}$ divides
$n$, we have that for every node $x$ of a curve $\cC$ underlying a
twisted $G_g^{(k),n}$-cover, the order of $\Gamma_x$   divides $n$.

The following 
 is a restatement of a result of Pikaart and De Jong 
(\cite{Pikaart-DeJong}, Theorem 3.1.3 part (4)).

\begin{proposition} Let $G=G_g^{(k),n}$ with $k\geq 4$ and $\gcd(n,6)=1$. Then 
the group homomorphism $M \to \Out(G)$ is injective.
\end{proposition}

As an immediate outcome we have

\begin{theorem}
Suppose  $k\geq 4$ and $\gcd(n,6)=1$. Then the moduli stack
$\TCtei{G_g^{(k),n}}$ is a smooth  
projective scheme over $\ZZ[1/n]$ admitting a flat finite morphism to
$\ocM_g$.  
\end{theorem}

We remark that, in this case, the space $\TCtei{G_g^{(k),n}}$ coincides with
the space ${}_G\obM_g$ (and the stack ${}_G\ocM_g$) of Pikaart and De Jong
(\cite{Pikaart-DeJong}, Section 2.3.5).

\subsection{A fine moduli space of $G$ covers}\label{Sec:fine-covers}

We are now ready to describe a particular type of finite group $G$ for which
the stack of connected admissible $G$-covers is representable; in
particular  the center of $G$ is trivial. The construction is closely
related to that of Looijenga \cite{Looijenga}, and we rely on some of
his arguments in our proofs.

We start with some auxiliary constructions. 

Let $p_1,p_2$ be two distinct primes. Fix a smooth complex curve $C$ of genus
$g>1$. Let $C_{p_1} \to C$ be the maximal finite abelian \'etale cover of
exponent $p_1$. 
The Galois group $G_{p_1}$ of $C_{p_1}$ over $C$ is $H_1(C, 
\ZZ/p_1\ZZ)$, which is isomorphic to $(\ZZ/p_1\ZZ)^{2g}$. 
The genus $g_{p_1}$ of $C_{p_1}$ is $p_1^{2g}(g-1)+1$.
Let $C_{\prym{p_2}{p_1}}\to C_{p_1}$ be the maximal finite abelian
\'etale cover of exponent $p_2$. The Galois group $H_{\prym{p_2}{p_1}}$ of
$C_{\prym{p_2}{p_1}}\to C_{p_1}$ is $H_1(C_{p_1}, 
\ZZ/p_2\ZZ)$, which is isomorphic to $(\ZZ/p_2\ZZ)^{2g_{p_1}}$.
The Galois group $G_{\prym{p_2}{p_1}}$ of the characteristic cover
$C_{\prym{p_2}{p_1}}\to C$ 
sits in an exact sequence:
$$ 1 \to H_{\prym{p_2}{p_1}} \to G_{\prym{p_2}{p_1}} \to
G_{p_1} \to 1.$$ 
Since $p_1 \neq p_2$, this is in fact a semidirect product. Moreover,
there is a canonical $G_{p_1}$-equivariant direct sum decomposition 
$$H_{\prym{p_2}{p_1}} = H_{\prym{p_2}{p_1}}^i\oplus
 H_{\prym{p_2}{p_1}}^n,$$
where the first factor consists of the $G_{p_1}$-invariants, and the
second of non-invariants. Explicitly, for $h\in H_{\prym{p_2}{p_1}}$,
we have the decomposition $h = h^i + h^n$ where $h^i =
\frac{1}{p_1^{2g}} \sum_{g\in G_{p_1}} h^g$. The subgroup
$H_{\prym{p_2}{p_1}}^i$ clearly lies in the center of
$G_{\prym{p_2}{p_1}}$. Denote the quotient by
$G_{\prym{p_2}{p_1}}^n$, and the corresponding covering 
$$C_{\prym{p_2}{p_1}}^n \to C.$$

\begin{lemma}\label{Lem:group-is-perfect}
\begin{enumerate} 
\item The action of $G_{p_1}$ on $H_{\prym{p_2}{p_1}}$ is effective.
\item The center of  $G_{\prym{p_2}{p_1}}$ is
$H_{\prym{p_2}{p_1}}^i$, and the center of $G_{\prym{p_2}{p_1}}^n$ is
trivial. 
\item The group $H_{\prym{p_2}{p_1}}^i$ is canonically isomorphic to
$G_{p_2}$. 
\item We have a diagram with cartesian squares:
$$\begin{array}{ccccc}
 C_{\prym{p_2}{p_1}} & \to&C_{p_1} \mathbin{\mathop\times\limits_C}
C_{p_2}&\to&C_{p_2} \\ 
\dar 			&  &\dar		&     &\dar    \\
C_{\prym{p_2}{p_1}}^n&\to &C_{p_1}		&\to  & C 
\end{array}$$
\end{enumerate}
\end{lemma}

\proof Clearly $H_{\prym{p_2}{p_1}}^n \subset G_{\prym{p_2}{p_1}}$
is a normal subgroup, and the quotient is abelian with exponent
$p_1p_2$. Thus its $p_2$-Sylow subgroup $H_{\prym{p_2}{p_1}}^i$
has order $\leq p_2^{2g}$. This implies that
$H_{\prym{p_2}{p_1}}^n$ is nontrivial. In particular this means that
some element of $G_{p_1}$ acts nontrivially on $H_{\prym{p_2}{p_1}}$.

Since the covering is characteristic, the  automorphisms group of
the fundamental group of $C$ acts on $G_{\prym{p_2}{p_1}} \to
G_{p_1}$. This automorphism group contains the Teichm\"uller group,
and  it is well known\marg{Precise citation needed!} that the
image of this group in $\Aut G_{p_1}$ is the symplectic
group, which is transitive on nonzero elements. It follows that {\em
every} nonzero element of 
$G_{p_1}$ acts nontrivially on $H_{\prym{p_2}{p_1}}$, giving (1). It
also follows that an element of the center maps to  $0\in G_{p_1}$,
therefore it is in $H_{\prym{p_2}{p_1}}$, and since it commutes with
$G_{p_1}$ 
it is in $H_{\prym{p_2}{p_1}}^i$, giving (2).

The curve $C_{p_1} \times_CC_{p_2}$ is an abelian  connected
$G_{p_2}$-covering of $C_{p_1}$, which means that $H_{\prym{p_2}{p_1}}^i$
has order $\geq p_2^{2g}$. Combined with the previous inequality we
get equality, giving (3) and (4). \qed

By switching the roles of $p_1$ and $p_2$ we can extend the fiber
diagram as follows:
$$\begin{array}{ccccc}
 C_{\{p_1,p_2\}} & \to & C_{\prym{p_1}{p_2}}&\to&
C_{\prym{p_1}{p_2}}^n\\   
\dar 			&  &\dar		&     &\dar    \\
 C_{\prym{p_2}{p_1}} & \to&C_{p_1} \mathbin{\mathop\times\limits_C}
C_{p_2}&\to&C_{p_2} \\ 
\dar 			&  &\dar		&     &\dar    \\
C_{\prym{p_2}{p_1}}^n&\to &C_{p_1}		&\to  & C 
\end{array}$$
The resulting curve $C_{\{p_1,p_2\}}$ is connected since the degrees
in the top left square are relatively prime. We denote the Galois
group of the covering $C_{\{p_1,p_2\}} \to C$ by 
$G_{\{p_1,p_2\}}$. It is clearly a characteristic quotient of the
fundamental group, isomorphic to
$G_{\prym{p_2}{p_1}}^n\times G_{\prym{p_1}{p_2}}^n$. Its exponent is
$p_1p_2$. 

We use the following Proposition:

\begin{proposition}\label{Prp:two-non-sep} Let $C$ be a stable curve
of genus $g>1$ and let $ 
C_{p_1}\to C$ be a 
connected admissible $G_{p_1}$-cover. Then no two nodes of $C_{p_1}$
separate the curve.
\end{proposition}

\proof The argument of Looijenga in \cite{Looijenga}, Proposition 2,
works word for word. \qed

\begin{lemma}\label{Lem:index-p1p2} Let $C$ be a stable curve of genus
$g>1$ and let $ 
C_{\{p_1,p_2\}}\to C$ be a 
connected admissible $G_{\{p_1,p_2\}}$-cover. Then at every node of
$C$, the covering has index $p_1p_2$. 
\end{lemma}

\proof The Proposition 
shows in particular that every node of $C_{p_1}$ is nonseparating. Our
Proposition 
\ref{Prop:pre-level} implies that the cover $C_{\prym{p_2}{p_1}}
\to C_{p_1}$ has 
index $p_2$. Similarly $C_{\prym{p_1}{p_2}} \to C_{p_2}$ has index
$p_1$. This implies that $p_1p_2$ divides the index of
$C_{\{p_1,p_2\}}\to C$.
On the other hand the index divides the exponent $p_1p_2$ of
$G_{\{p_1,p_2\}}$, giving equality. \qed

\begin{lemma}\label{Lem:action-p1p2}  Let $C$ be a stable curve of
genus $g>1$ and let $ 
C_{\{p_1,p_2\}}\to C$ be a 
connected admissible $G_{\{p_1,p_2\}}$-cover. Consider the homomorphism
 $\partial:\prod_{e \in C\sing} \ZZ \to \Out (G_{\{p_1,p_2\}})$ via Dehn
twists. If the image of an element $u\in \prod_{e \in C\sing}
\ZZ$ in $\Out (G_{\{p_1,p_2\}})$ is trivial, then $u\in
p_1p_2\prod_{e \in C\sing} \ZZ$. 
\end{lemma} 

\proof Denote by $E_0$ the set of separating nodes of $C$ and by $E_1$
the nonseparating nodes. We now follow the notation of
\cite{Looijenga}: we write $e\in C\sing$ for a node of
$C$ (instead of $x$ used earlier); for such a node we denote the corresponding
component of $u$ by $u_e$; a node of $C_{p_1}$  is
denoted $\tilde e$, and we indicate the condition that it lie over $e$
by $\tilde e/e$. We denote by $[\tilde e]$ the class of the
corresponding vanishing cycle in  $H_1(C_{p_1}^*, \ZZ)$ (with some  choice
of orientation), where
$C_{p_1}^*$ is a nearby smooth fiber in a deformation of $C_{p_1}$.

According to Looijenga's discussion in
\cite{Looijenga}, 
p. 287-288 the action of the element $p_1 u$ on $v\in H_1(C_{p_1}^*, \ZZ)$ is
$$v \mapsto v+ \sum_{e\in E_0} p_1 u_e \sum_{\tilde e/e} (v,[\tilde e])
[\tilde e]  +  \sum_{e\in E_1} u_e \sum_{\tilde e/e} (v,[\tilde e])
[\tilde e].$$  If  $\partial u$  is
trivial, then $u$ acts on $H_1(C_{p_1}^*, \ZZ/p_2\ZZ)$ as an element of
$G_{p_1}$, therefore $p_1u$ acts trivially on $H_1(C_{p_1}^*, \ZZ/p_2\ZZ)$. 
Proposition \ref{Prp:two-non-sep} says that Baclawski's Lemma
(\cite{Looijenga}, p. 286) applies, therefore $p_2 | p_1 u_e$ for
$e\in E_0$ and $p_2 | u_e$ for $e\in E_1$. Since $p_2\neq p_1$ it
follows that $p_2|u_e$ 
for all $e$, i.e. $u\in   p_2\prod_{e \in C\sing} \ZZ$. Reversing the
roles of $p_1$ and $p_2$ we get the Lemma.\qed

\begin{theorem} The automorphism group of any connected admisible
$G_{\{p_1,p_2\}}$-cover is trivial.
\end{theorem}

\proof By Lemma \ref{Lem:char-dehn} we have a surjection $\delta:\prod_{e
\in C\sing} \ZZ\to M$ compatible with the 
action by outer automorphisms, in other words, the homomorphism
$\partial:\prod_{e \in C\sing} 
\ZZ \to \Out (G_{\{p_1,p_2\}})$ factors through $\delta$. Thus  
$$\Ker\ \delta \ \subset \ \Ker\ \partial.$$
By Lemma 
\ref{Lem:index-p1p2} we have $$\Ker\ \delta \ = \   p_1p_2\prod_{e \in
C\sing} \ZZ.$$ By lemma \ref{Lem:action-p1p2} we have that 
 $$\Ker\ \partial \ \subset \ p_1p_2\prod_{e \in
C\sing} \ZZ.$$
 Combining the statements we get equality.  Thus the
automorphism group of the cover is the center of the group, which by
Lemma \ref{Lem:group-is-perfect} is
trivial.\qed

\begin{remark} We note that similar results can be obtained for covers
of pointed curves, for instance using the reduction methods of
\cite{Boggi-Pikaart}. 
\end{remark}
\marg{In final version we should add  the example of  a connected
$G$-cover with $Z(G)$ trivial and $A^G$ nontrivial. Maybe also
discuss infinitely twisted curves.}

\appendix\section{Some remarks on \'etale cohomology of
Deligne--Mumford stacks}

When we refer to a sheaf on a stack or algebraic space $\cM$, we will
always mean a sheaf in the small \'etale site of $\cM$, whose objects are
\'etale morphism locally of finite type $U \to \cM$, where $U$ is a scheme.

\begin{proposition}\label{Prop:hdi-to-moduli} Let $\cM$ be a separated tame
finitely presented Deligne--Mumford stack over a scheme, with moduli space
$\pi : \cM \to M$. Let $p : \Spec\Omega \to \cM$ be a geometric point of
$\cM$, $\Gamma$ its stabilizer, $q = \pi \circ p : \Spec\Omega \to M$ its
image in $M$. Let $F$ be a sheaf on $\cM$; then there is a canonical
isomorphism of groups between the stalk $(\bR^i\pi_*F)_q$ of the
$i^{\operatorname{th}}$ higher direct image sheaf of $F$ at $p$ with the
$i^{\operatorname{th}}$ cohomology group $H^i(\Gamma, F_p)$.
\end{proposition}

\proof This statement is local in the \'etale topology on $M$, so we can 
make
an \'etale base change and assume that $\cM$ is a quotient $[U/\Gamma]$, where
$U$ is a connected scheme and $p$ is a geometric point which is fixed by
$\Gamma$; then the sheaf $F$ is a $\Gamma$-equivariant sheaf on $U$, and 
$M = U/\Gamma$. Let $\Spec\Omega \to V \to U/\Gamma$ be an \'etale neighborhood
of $q$ in $U/\Gamma$; then $[V\times_{U/\Gamma} U/\Gamma] = V\times_{U/\Gamma}
[U/\Gamma]$. There is a spectral sequence
$$
E_2^{ij} = H^i\bigl(\Gamma, H^j(V\times_{U/\Gamma} U, F)\bigr)
\Longrightarrow
H^{i+j}(V\times_{U/\Gamma} [U/\Gamma], F)
$$
which is the \v Cech-to-global cohomology spectral sequence for the covering $U
\to [U/\Gamma]$, as in \cite{Milne}, Proposition~2.7, with the same proof. Now
let us go to the limit over all \'etale neighborhoods $\Spec\Omega \to V \to
U/\Gamma$; since $\Gamma$ is a finite group, its cohomology groups commute with
direct limits, so we have
$$
\varinjlim_V H^i\bigl(\Gamma, H^j(V\times_{U/\Gamma} U, F)\bigr) =
H^i\bigl(\Gamma, \varinjlim_V H^j(V\times_{U/\Gamma} U, F)\bigr).
$$

But $ \varinjlim_V H^j(V\times_{U/\Gamma} U, F)$ is 0 for $j > 0$ and is
$F_p$ for $j = 0$, while the limit of the abutment of the spectral sequence
is precisely $(\bR^i\pi_*F)_q$. \qed

Now we prove the proper base change theorem for tame Deligne--Mumford stacks.

A sheaf $F$ on a stack $\cM$ is {\em torsion} if for any \'etale
morphism $U
\to \cM$ from a quasicompact scheme $U$ the group $F(U)$ is torsion.

\begin{theorem}\label{Th:proper-base-change}
Let $f : \cM \to S$ be a proper morphism from a tame Deligne--Mumford
stacks to a scheme, and let
$$
\begin{array}{ccc}
\cM'    & \stackrel{\psi}{\lrar} & \cM    \\
\dar f' &                        & \dar f \\
S'      & \stackrel{\phi}{\lrar} & S 
\end{array}
$$
be a cartesian diagram.
Let $F$ be a torsion sheaf on $\cM$. Then the natural base change  
homorphism of sheaves $\phi^*\bR^if_* F \to \bR^if'_*\psi^* F$ is an
isomorphism.
\end{theorem}

\proof When $\cM$ is a scheme, this is the usual proper base change theorem for
\'etale cohomology, as in \cite{Milne}, Corollary~2.3. This also holds when
$\cM$ is an algebraic space, with the same proof. We will reduce the general
case to the case of algebraic spaces.

First of all, if $S$ is the moduli space of $\cM$, the statement follows
easily from \ref{Prop:hdi-to-moduli}. In general, factoring through the
moduli spaces we get a cartesian diagram
$$
\begin{array}{ccc}
\cM'      & \stackrel{\psi}{\lrar}  & \cM      \\
\dar \pi' &                         & \dar \pi \\
M'        & \stackrel{\rho}{\lrar}  & M        \\
\dar g  ' &                         & \dar g   \\
S'        & \stackrel{\phi}{\lrar}  & S        \\
\end{array}
$$
such that $g \circ \pi = f$ and $g' \circ \pi' = f'$. The base change
formula holds for $pi$ and $g$, so we have
$$
\phi^* \bR^ig_* \bR^j\pi_* F =  \bR^ig'_* \rho^* \bR^j\pi'_* F =
 \bR^ig'_* \bR^j\pi'_* \psi^* F.
$$

We also have a morphism of $E_2$ spectral sequences 
$$\begin{array}{ccc}
\phi^* \bR^ig_* \bR^j\pi_* F  & \Longrightarrow & \phi^*\bR^{i+j}f_*F  \\
\dar                          &                 & \dar                 \\
\bR^ig'_* \bR^j\pi'_* \phi^*F & \Longrightarrow & \bR^{i+j}f'_*\psi^*F
\end{array}$$
where the columns are base change maps. Since the left hand column is an
isomorphism, so is the right hand column. \qed

\end{document}